\documentclass[absolute]{mymathart}
\usepackage[theorems]{mymathmacros}
\usepackage[usenames,dvipsnames]{color}
\usepackage{subfig}
\usepackage[shortlabels]{enumitem}
\usepackage{hyperref}

\title[Hyperbolic entire functions with bounded Fatou components]{Hyperbolic entire functions with bounded Fatou components}
\author{Walter Bergweiler}
\address{Mathematisches Seminar, Christian-Albrechts-Universit\"{a}t zu Kiel, 
24098 Kiel, Germany}
\email{bergweiler@math.uni-kiel.de}
\author{N\'uria Fagella}
\address{Dep.~de Matem\`atica Aplicada i An\`alisi, Universitat de Barcelona, Gran Via 585, 08007 Barcelona, Spain}
\email{fagella@maia.ub.es}
\author{Lasse Rempe-Gillen}
\address{Dept.\ of Mathematical Sciences, University of Liverpool, Liverpool L69 7ZL, UK}
\email{l.rempe@liverpool.ac.uk}
\thanks{The second author was supported by 
Polish NCN grant decision DEC-2012/06/M/ST1/00168 as well as grants
2009SGR-792 and MTM2011-26995-C02-02. The third author was supported by a Philip Leverhulme Prize.}
\date{\today}

\newcommand{\LP}{\text{\it LP}}

\newtheorem{obs}[thm]{Observation}

\newcommand{\interior}{\operatorname{int}}

\numberwithin{equation}{section}

\renewcommand{\P}{\mathcal{P}}

\newcommand{\B}{\mathcal{B}}

\renewcommand{\H}{\mathbb{H}}

\renewcommand{\Ch}{\widehat{\C}}

\newcommand{\tn}{\tiny}

\begin{document} 
 \begin{abstract}
  We show that an invariant Fatou component of a 
   hyperbolic transcendental entire function is a 
   Jordan domain 
   (in fact, a quasidisc)
    if and only if it contains only finitely many critical points and no asymptotic curves. 
    We use this theorem to prove criteria for the boundedness of Fatou components and
    local connectivity of Julia sets for hyperbolic entire functions,
    and give examples that demonstrate that our results are optimal.
    A particularly strong dichotomy is obtained in the case of a function with 
    precisely two critical values.
 \end{abstract}
 
 \maketitle

 \section{Introduction}
 
   Dynamical systems that are \emph{hyperbolic} (or ``Axiom A'' in Smale's terminology)
   exhibit, in a certain sense, the simplest possible behaviour. (For the formal definition
   of hyperbolicity in our context, see Definition \ref{defn:hyperbolicity} below.) In any given setting,
   understanding hyperbolic systems is the first step on the way to studying more
   general types of behaviour. Furthermore, in many one-dimensional situations, hyperbolic
   dynamics is either known or believed to be topologically generic (see e.g.\ \cite{lyubichdensity,graczykswiatek,KSS2,RGvS}), and hence
   many systems are indeed hyperbolic.

  In the iteration of complex polynomials $p\colon\C\to\C$, the dynamics of 
   hyperbolic functions has been essentially completely understood since the 
   seminal work of Douady, Hubbard and Thurston in the 1980s. In particular, these
   can be classified~--~in a variety of ways~-- using finite combinatorial objects such as
    ``Hubbard trees''. Typically, any qualitative question about the iterative behaviour
   of the map in question
    can be answered from this encoding. 

 In addition to polynomial and rational iteration, the dynamical study of 
   transcendental entire functions (i.e.non-polynomial holomorphic
   self-maps of the complex plane) is currently receiving increased interest, partly
    due to intriguing connections with deep aspects of the polynomial theory.
    (We refer to the introduction of \cite{strahlen} for a short discussion.) 
    However, until recently there were only a small number of specific cases where 
   hyperbolic behaviour had been understood in detail (cf.\ \cite{aartsoversteegen,bhattacharjee,schleicherzimmerperiodic} and 
    \cite[Corollary 9.3]{topescaping}).

  Indeed, it turns out that, even restricted to the hyperbolic case, entire functions can be incredibly diverse: 
    for example, while for many such maps, the Julia sets are known to contain curves along which 
    the iterates tend to infinity, there are also (hyperbolic) examples where this is not the case 
    \cite{strahlen}. Similarly, for some hyperbolic maps there are natural conformal  
    measures in the sense of Sullivan, with associated invariant measures 
    \cite{mayerurbanski}, while for others such measures cannot exist \cite{fullhypdimension}. 
    Nonetheless, it was proved recently \cite[Theorems~1.4 and~5.2]{rigidity}
    that, in any given family of entire functions, the behaviour of hyperbolic functions can essentially be
    described completely, in terms of a certain topological model (which however 
   depends on the family in question).

 A disadvantage of this description is that it is not very explicit. To explain what we mean by this, 
   and to introduce 
   the main question treated in our article, 
   we first provide some of the definitions that were deferred above. 
    A point $s\in\C$ is called a \emph{singularity of the inverse function}
    $f^{-1}$ if $s$ is either a \emph{critical value} (the image of a critical point) or
    a finite \emph{asymptotic value}. The latter means that there is a path $\gamma$ to infinity whose image
    ends at $s$; the curve $\gamma$ is then called an \emph{asymptotic curve}. The set of such singularities is 
    denoted by $\sing(f^{-1})$. Then we call 
      \begin{equation}\label{eqn:classB}
      \B := \{f\colon\C\to\C \text{ transcendental entire}\colon \sing(f^{-1})\text{ is bounded}\} 
     \end{equation}
  the \emph{Eremenko-Lyubich class} (compare \cite{eremenko_lyubich_2}). 

 \begin{defn}[Hyperbolicity]\label{defn:hyperbolicity}
  A transcendental entire function $f\colon\C\to\C$ is called \emph{hyperbolic} if $f\in\B$ and 
   furthermore every element of $S(f) := \overline{\sing(f^{-1})}$ belongs to the basin of
   some attracting periodic cycle of $f$. 
 \end{defn}
Equivalently, $f$ is hyperbolic if and only if the \emph{postsingular set}
  \[ \P(f) := \overline{\bigcup_{j\geq 0} f^j(\sing(f^{-1}))} \]
  is a compact subset of the Fatou set (see Proposition \ref{prop:hyproperties}). 
  Recall that the \emph{Fatou set}, $F(f)$, consists of those points whose behaviour under iteration is stable;
    more precisely, it contains exactly those points where the family of iterates 
    $(f^n)_{n\geq 0}$ 
  is equicontinuous with respect to the spherical metric. 
   The complement $J(f) = \C\setminus F(f)$ is
  called the Julia set; here the dynamics is unstable.

\begin{figure}[htb!]
  \subfloat[$f(z)=e^z+2 + \frac{\pi i}{2}$\label{subfig:expo}]
       {\includegraphics[width=.39\textwidth]{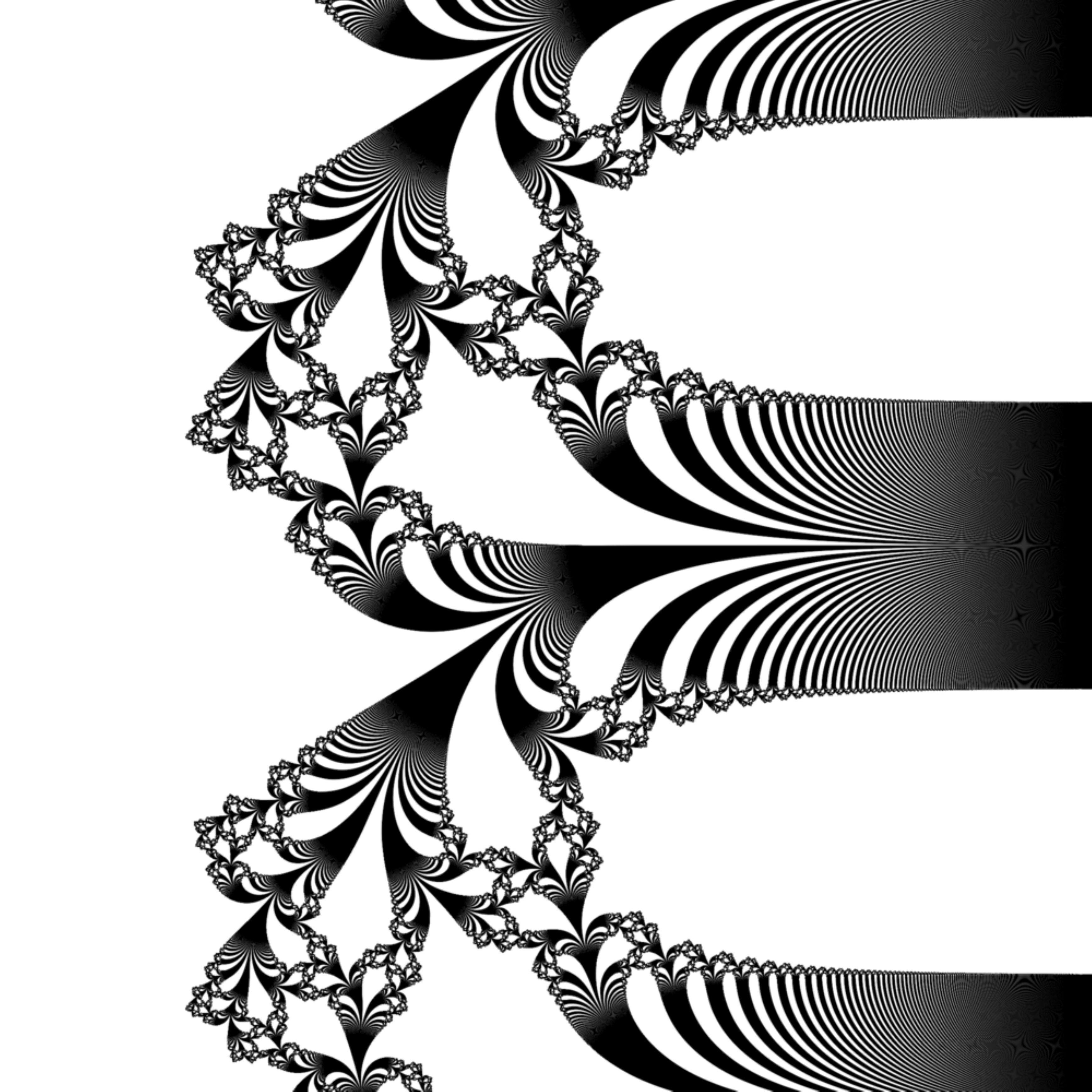}}
    \hfill
  \subfloat[$f(z)=-z^2\exp(1-z^2)$]{\includegraphics[width=.585\textwidth]{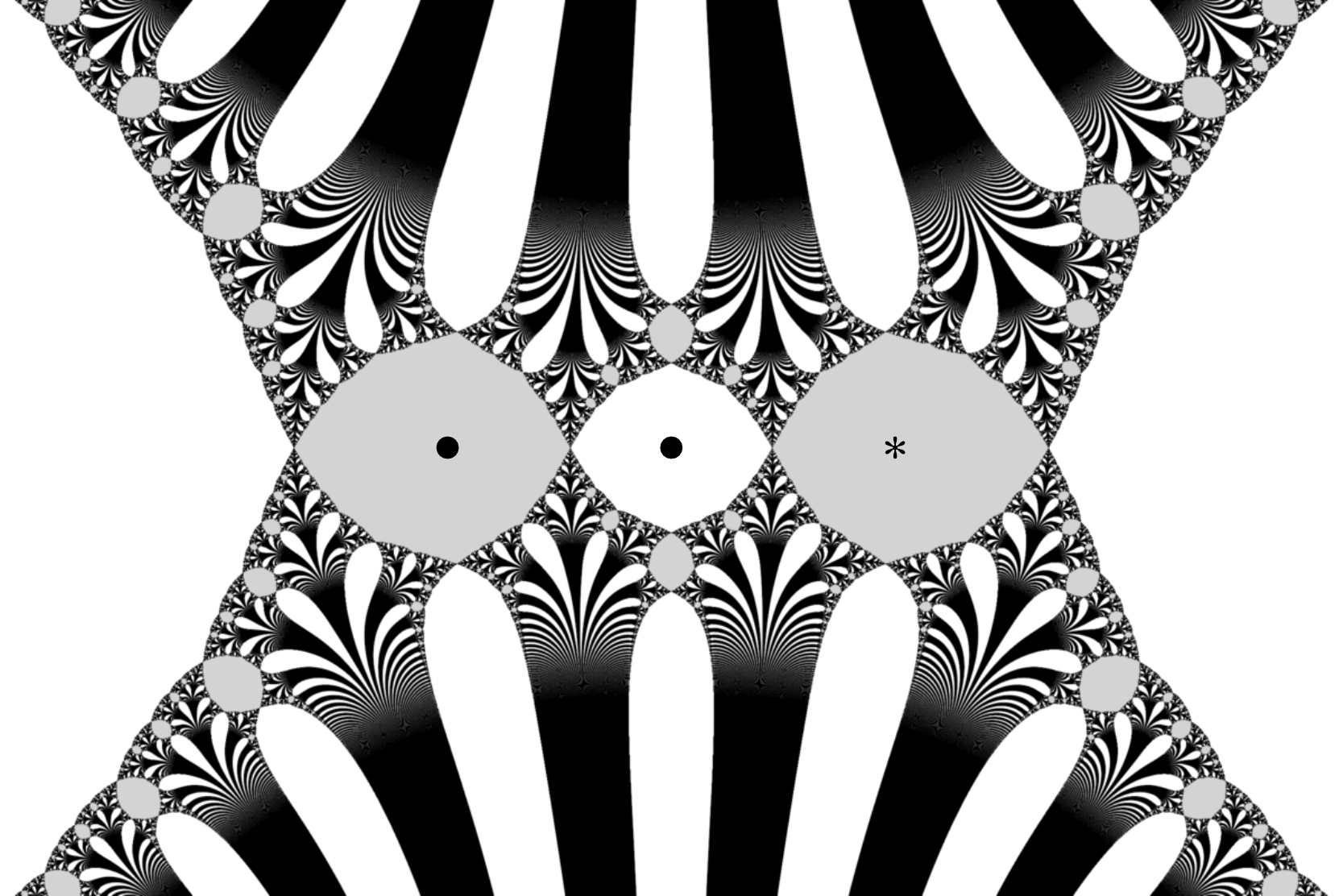}\label{subfig:maclanevinberglimit}}
 \caption{\label{fig:asymptotic} \small Two entire functions with asymptotic values
   and unbounded Fatou components (the Julia set is drawn in black; Fatou components are in grey and white). On the left is an exponential map;
   the Fatou set consists of the basin of an attracting orbit of period $3$. Every Fatou
   component $U$ is unbounded and $\partial U$ is not locally connected. On the right
   is a function that plays a crucial role in our construction of Example \ref{example1}. 
  Here there
   are superattracting fixed points at $0$ and $-1$ (marked with filled circles), 
   and $0$ is an asymptotic value. The basin of $0$ is coloured white and the basin of $-1$
   is coloured grey~--
   every Fatou component is a Jordan domain, but all pre-periodic components of the 
   basin of $0$ are unbounded,
    and the Julia set is not locally connected. Here and in subsequent images, 
   non-periodic critical points (in this case, the point $1$) are marked by asterisks.}
\end{figure}

\begin{figure}[htb!]
  \subfloat[\label{subfig:cosdisjoint}$F(f)$ connected]{
  \includegraphics[width=.49\textwidth]{%
   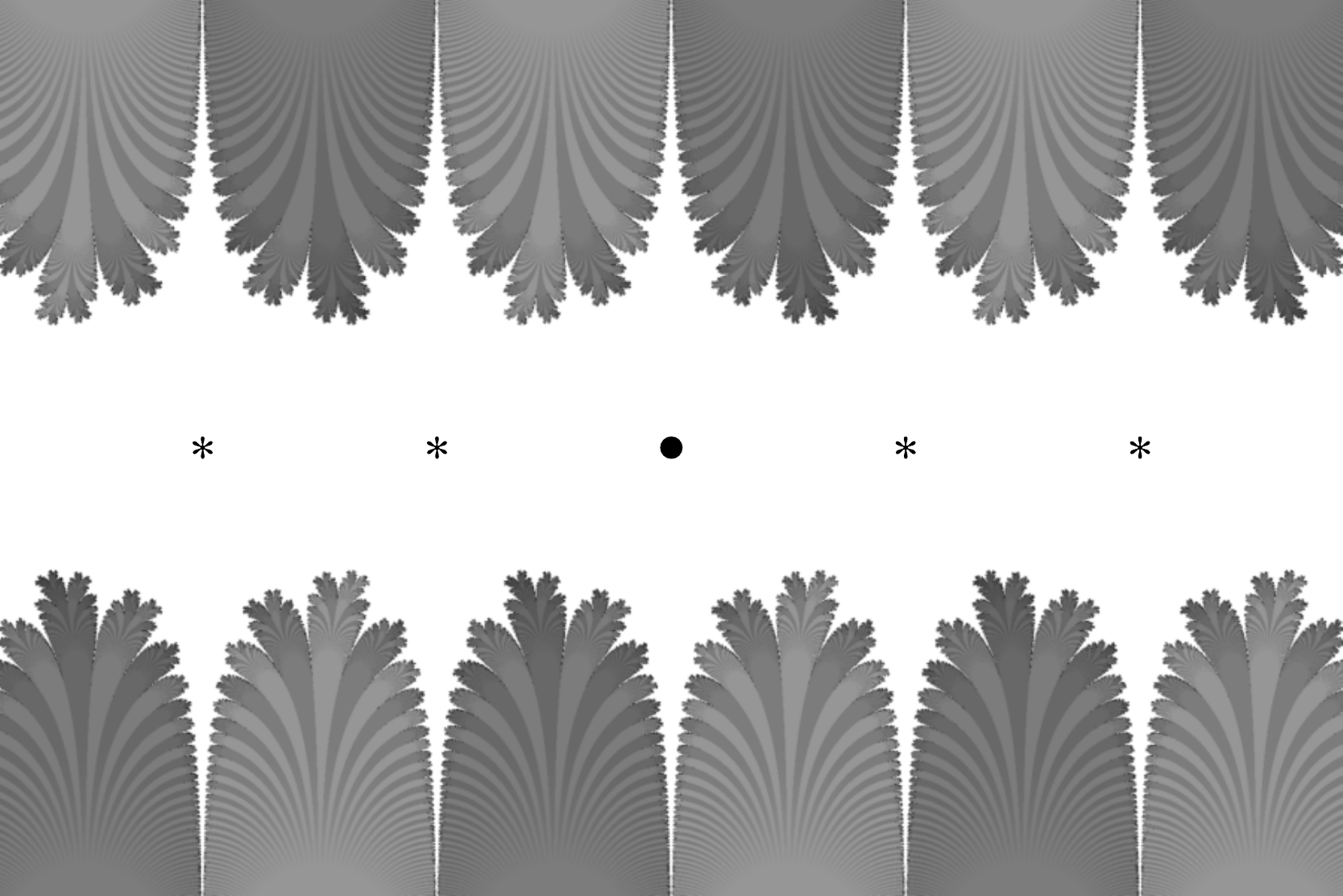}}\hfill
  \subfloat[\label{subfig:cosunbounded}Unbounded components]{%
  \includegraphics[width=.49\textwidth]{%
   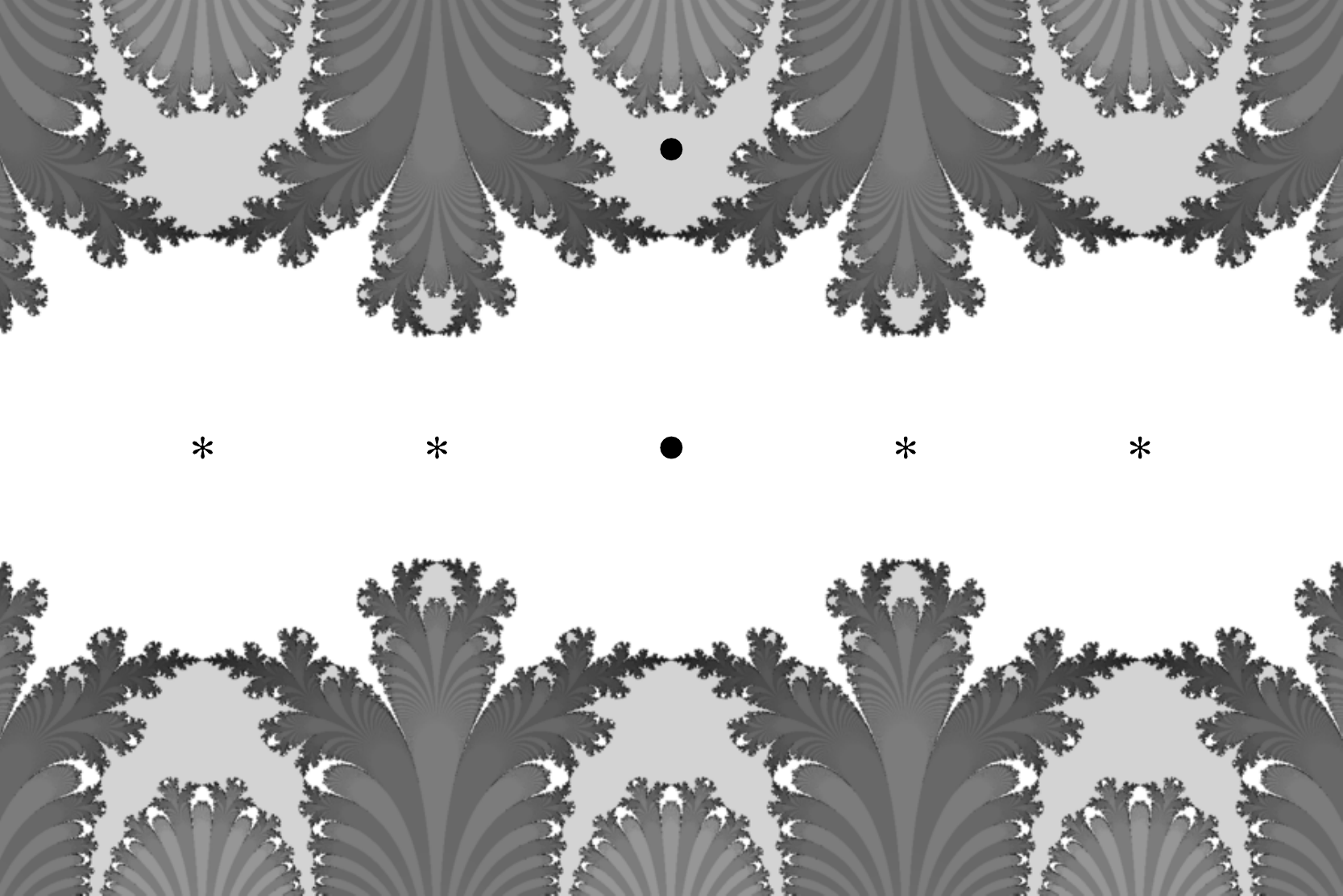}}\\
  \subfloat[\label{subfig:cospreperiod}Locally connected Julia set]{%
  \includegraphics[width=.49\textwidth]{%
    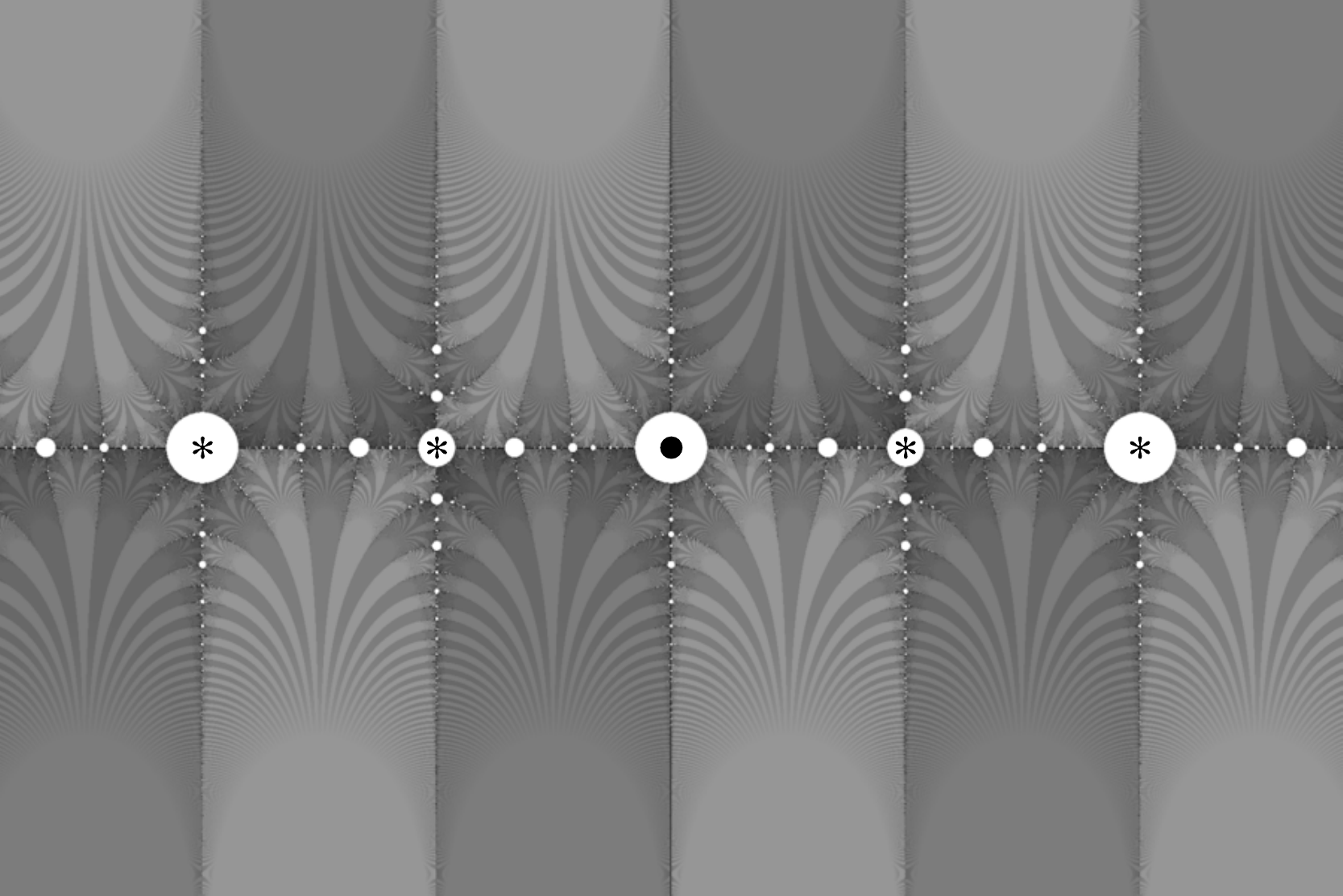}}\hfill
  \subfloat[\label{subfig:cosbasilica}Another locally connected Julia set]{%
  \includegraphics[width=.49\textwidth]{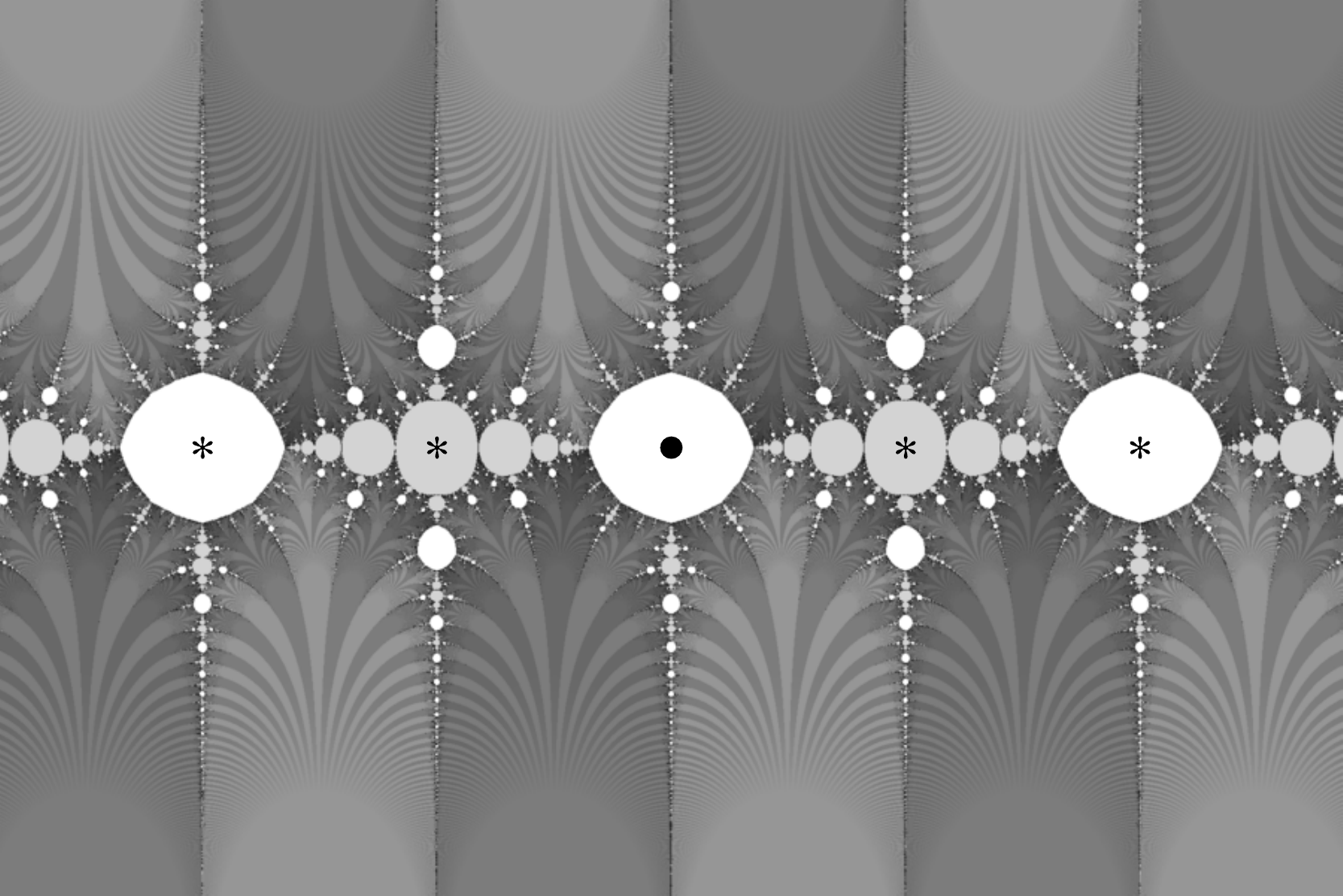}}
  \caption[Julia sets i the cosine family]{%
   \label{fig:sine} \small Julia sets (in dark grey) of some maps in the cosine family,
     $z\mapsto a\cos(z)+b$, illustrating Theorem \ref{thm:jordantwo} 
     and Corollary
     \ref{cor:lctwo}.{\protect\footnotemark}
 (Note that these can easily be reparametrised as 
     $z\mapsto \sin(a'z + b')$ for suitable choices of $a',b'$.)
    The two maps in the top row have unbounded Fatou components
     and their Julia sets are not locally connected. The maps on the bottom
      line both have locally connected Julia sets.}
\end{figure}

  If a hyperbolic entire function $f$ has an asymptotic value, then it follows immediately that
   some (but not necessarily all) connected components of $F(f)$ are unbounded, see 
   Figures \ref{fig:asymptotic}\subref{subfig:expo} and \ref{fig:asymptotic}\subref{subfig:maclanevinberglimit}.  
   On the other hand, there are also examples where 
   all Fatou components appear to be bounded Jordan domains; compare Figures \ref{fig:sine}\subref{subfig:cospreperiod} and
           \ref{fig:sine}\subref{subfig:cosbasilica}.
 Unfortunately, the above-mentioned description
  from \cite{rigidity} does not allow us to determine when this is the case, and hence the problem remained open 
   even for 
   rather simple explicit cases such as the \emph{cosine family}, see below. The following result gives a complete answer to this 
   question, and hence provides another step towards the understanding of hyperbolic transcendental entire
   dynamics. 

  \begin{thm}[Bounded Fatou components] \label{thm:alljordan}
    Let $f\in\B$ be hyperbolic. Then the following are equivalent:
     \begin{enumerate}[{\normalfont(1)}]
       \item every component of $F(f)$ is bounded;\label{item:allbounded}
       \item $f$ has no asymptotic values and every component of $F(f)$ contains at most
          finitely many critical points.\label{item:allcritbounded}
     \end{enumerate}
  \end{thm}
\begin{remark}
 If either (and hence both) of these conditions are true, then in fact all Fatou components
  are bounded quasidiscs, see Corollary \ref{cor:quasidiscs}.
  (A \emph{quasidisc} is a Jordan domain that is the image of the open unit disc under some quasiconformal  homeomorphism of the Riemann sphere.)
\end{remark}

Condition~\ref{item:allcritbounded} of Theorem~\ref{thm:alljordan} can usually be verified in a straightforward manner for specific hyperbolic functions. 
This is especially easy when $\sing(f^{-1})$ is finite, $f$ has no asymptotic values
 and we know that every Fatou component contains at most one critical value (for example, because different
 critical values converge to different attracting periodic orbits). We shall see 
  (in Proposition \ref{prop:covering}~\ref{item:onecriticalpoint}) 
  that in this case every component of $F(f)$ contains at most one critical point, and hence we obtain the following corollary:
\begin{cor}[One critical value per component]
\label{cor:onlyonecriticalpoint}
 Let $f$ be a hyperbolic entire function without asymptotic values. If every component of $F(f)$
   contains at most one critical value, then every component of $F(f)$ is bounded.
\end{cor}

\subsection*{Maps with two critical values}

   Let us consider what happens when we restrict the size of $\sing(f^{-1})$ further.
   If $\# \sing(f^{-1})=1$, then $f$ must be conjugate
  to an exponential map $z\mapsto \lambda e^z$. This family has been thoroughly studied since the 1980s;
 compare e.g.\ \cite{BakerRippon84,Devaney84,eremenko_lyubich_2,RemSch09} and the references therein. Since exponential maps have
   an asymptotic value at $0$ and no other singular values, in the hyperbolic case there are only unbounded Fatou components (see 
    Figure~\ref{fig:asymptotic}\subref{subfig:expo}).
\footnotetext{%
The maps in Subfigures \subref{subfig:cosdisjoint}, \subref{subfig:cospreperiod} and
   \subref{subfig:cosbasilica} have $-a=b=\lambda$, with
    $\lambda=3/4$, $\lambda=4\pi/3$ and $\lambda=2$, respectively. 
    For \subref{subfig:cosdisjoint} and \subref{subfig:cospreperiod} the
   superattracting fixed point at $0$ is the only attracting cycle, while in 
   \subref{subfig:cosbasilica}, there is an additional cycle of period $2$ whose basin is
   shown in light grey. In \subref{subfig:cosunbounded}, $a=4i/(1-\cosh4)$ and 
   $b=4i-a$. There 
   is a unique superattracting cycle $0\mapsto 4i\mapsto 0$. Points in $F(f)$ are
   coloured white and light grey depending on whether they take an even or odd number
   of iterations 
   to reach the Fatou component containing $0$. Superattracting cycles are indicated by
   black filled circles; non-periodic critical points are marked with asterisks.}
 
  Cases
   where $\# \sing(f^{-1})=2$
  include for instance the 
   \emph{cosine (or sine) maps} $z\mapsto \sin(az+b)$, where $a,b\in\C$, $a\neq 0$, with
   critical values at $\pm 1$ but no asymptotic values, as well as 
     the family $z\mapsto a z e^z+b$  with one critical and one asymptotic
    value, among many others (compare \cite{deniz2013,nurialimitingstandard}). 
   However, the class of entire maps with two inverse function singularities is far more general than suggested by these simple examples.
Indeed, there exist uncountably many essentially different families of entire functions
  with no asymptotic values and exactly two critical values;
 the same is true for functions with two asymptotic values, or one critical and one asymptotic value. 
   By this we mean that there exists an uncountable collection $\mathcal{F}$ of functions
of this type, such that none of the functions in $\mathcal{F}$ can be obtained from another by 
pre- and post-composition with plane homeomorphisms. (The existence of such a collection follows from the classical theory of \emph{line complexes}; compare 
  \cite[Chapter~7]{goldbergostrovskii} and also Observation \ref{obs:uncountable} below.)

Drasin \cite{Drasin07} and Merenkov \cite{Merenkov08} have constructed
  maps of this class that have irregular and arbitrarily fast growth, respectively. 
  More recently, Bishop \cite{bishopfolding} has described a method for constructing functions with no
   asymptotic values, two
   critical values and only simple critical points, having essentially arbitrary prescribed behaviour near infinity.
   These functions can have dynamical properties that are very different from those of the simple examples mentioned above. 
    For example, Bishop  shows \cite[Section~18]{bishopfolding} that the above-mentioned example from 
    \cite[Theorem 1.1]{strahlen}, where the \emph{escaping set} 
    \[ I(f) := \{z\in\C \colon f^n(z)\to\infty \} \]
     contains no unbounded continuous curves, can be realised within this class. 

Theorem~\ref{thm:alljordan} allows us to formulate a striking
     dichotomy when $\# \sing(f^{-1})=2$:
  \begin{thm}[Dichotomy for functions with two critical values]\label{thm:jordantwo}
  Let $f\colon\C\to\C$ be a transcendental entire function without finite asymptotic values
    and exactly two critical values. Assume furthermore that $f$ is hyperbolic, i.e.~both critical values tend to
    attracting periodic orbits of $f$ under iteration. Then either
   \begin{enumerate}[{\normalfont(1)}]
       \item \label{item:completelyinvariant} every connected component $U$ of $F(f)$ is 
         unbounded, and $\partial U$
         is not locally connected at any finite point, or 
        \item \label{item:jordan} every connected component of $F(f)$ is a bounded quasidisc.
    \end{enumerate}
  In case~\ref{item:completelyinvariant}, all critical points of $f$ belong to a single periodic Fatou component.
\end{thm}
Here ``local connectivity'' of a set $K$ at a point $z$ means that there are arbitrarily small
   connected 
   neighbourhoods of $z$ in $K$; we do not require these neighbourhoods to be open.
   (Sometimes this property is instead referred to as ``connected im kleinen''; compare e.g.\ 
   \cite{localconnectivity}.) Of course, the boundary of a quasidisc is locally connected at every point. Hence,
   for hyperbolic maps with two critical values, there are two extremely contrasting possibilities for the shape of all Fatou components.

The theorem appears to be new even for hyperbolic maps in the cosine family, except
  in the one-dimensional slice $z\mapsto \sin(az)$, where due to symmetry there is effectively only one free critical value.
     Here our result implies that
     all Fatou components are bounded when $f$ is hyperbolic, except for $|a| < 1$; this 
    was already stated
    by Zhang \cite[p.\ 2, third paragraph]{zhangsinecharacterization}, who mentions that it can be proved
    using polynomial-like mappings. Some special instances can also implicitly be found already
    in \cite{dominguezsienra}. The same statement is established in \cite[Prop.~6.3]{domfag} for all maps with 
    $|\re a|\geq \pi/2$,  even without the assumption of hyperbolicity; it is conjectured there that this also holds whenever $|a|\geq 1$.

\begin{figure}
{\def\svgwidth{\textwidth}
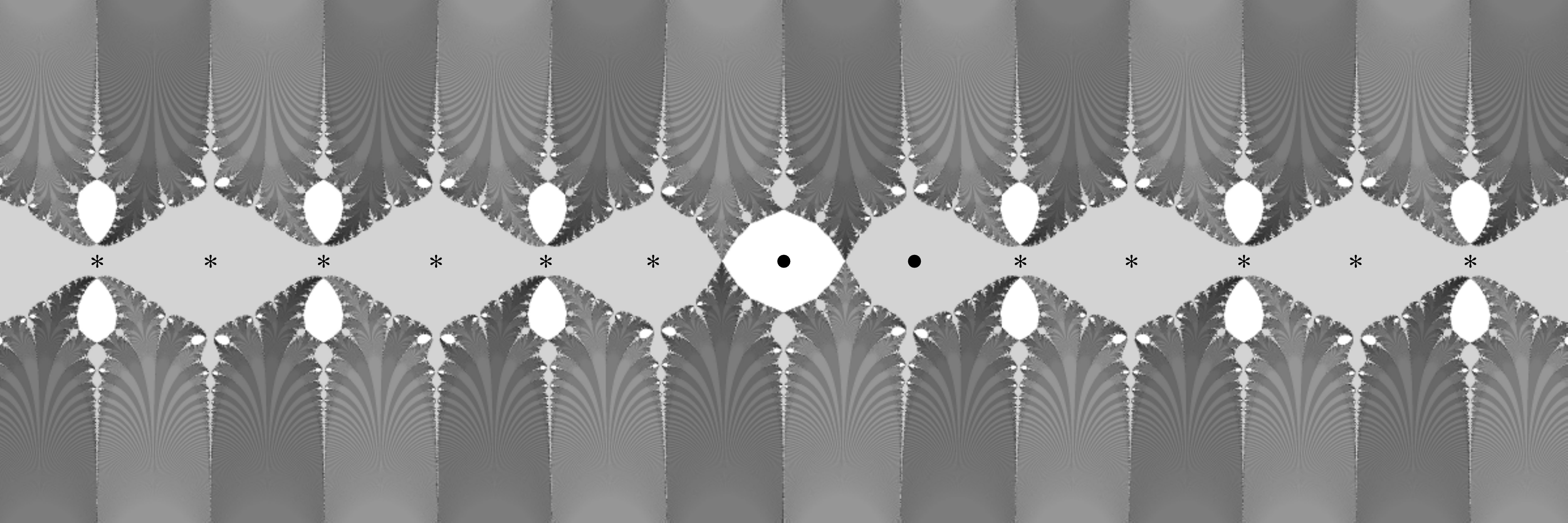}\\
{\def\svgwidth{\textwidth}
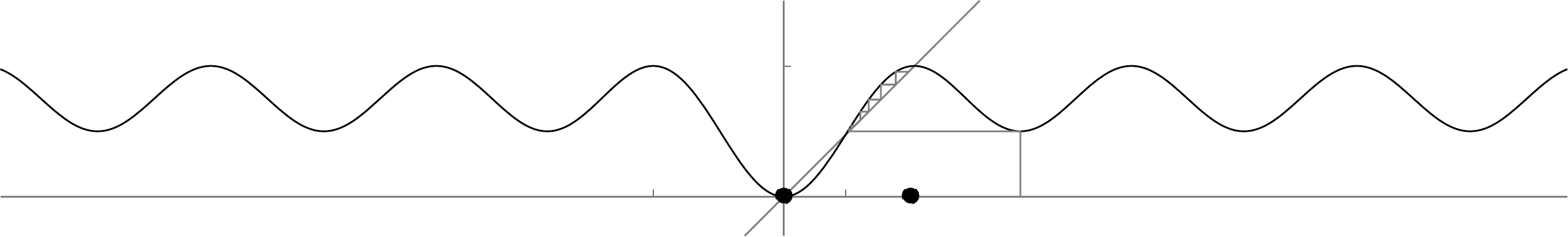}
\caption{\label{fig:cossqrt} \small The function from Example \ref{ex:cossqrt} (shown here for $u=3$) 
  has three critical values and two fixed Fatou components, one of which
    is a bounded Jordan domain, while the other is unbounded with non-locally-connected boundary. 
    The fixed points $0$ and $1$ are superattracting; the basin of $0$ is depicted in white, 
    the basin of $1$ in light grey, and the Julia set in darker tones of grey. 
    Non-periodic critical points are marked by asterisks.}
\end{figure}

Without further hypotheses,  Theorem \ref{thm:jordantwo} does not hold for larger numbers of critical values, as shown by the following example  (see Figure \ref{fig:cossqrt}).
  
 \begin{example}[Unbounded and bounded Fatou components]\label{ex:cossqrt}
 Define
 \[
 f(z)=
 \frac{u-\cos\sqrt{\left(\operatorname{arcosh}^2 u+\pi^2\right)z^2 -\operatorname{arcosh}^2 u}}{u+1}
 \] 
 where $u>1$. Then $f$ has no asymptotic values and three critical values $0$, $1$ and $c_u=(u-1)/(u+1)$.
 The points $0$ and  $1$ are superattracting fixed points and for large $u$ (in fact, for
 $u>2.7981186\dots$) the critical value $c_u$ and all positive critical points
 are contained in the immediate attracting basin of~$1$.
 The immediate basin of $0$ is a bounded quasidisc while the immediate basin of $1$ is unbounded and not
  locally connected at any finite point; see Figure \ref{fig:cossqrt}.
  \end{example}
 
 Even when the periodic Fatou components are Jordan domains, it is possible for some preimage components to be unbounded:
 \begin{example}[Unbounded preimages of bounded Fatou components]\label{ex:unbounded}%
 \label{example1}
There exists an entire function $f$, with three critical values $w_1,w_2,w_3$
and no asymptotic values, such that the following hold:
  \begin{enumerate}[{\normalfont(a)}]
   \item
   $w_1$ and $w_2$ are superattracting fixed points with immediate
basins bounded by a Jordan curve;
   \item
   $w_3$ is contained in the immediate basin of $w_1$;
   \item
   the preimage of the immediate basin of $w_1$ has an unbounded component.
  \end{enumerate}  
\end{example}

\subsection*{Local connectivity of Julia sets} 
  A key question in polynomial dynamics is whether a given Julia set is \emph{locally connected}.
   This is known to hold for large classes of examples,
    including all hyperbolic maps, 
    and implies a complete description of the topological dynamics  of the map in question
    (compare \cite{douadypincheddisk}).

 Local connectivity of Julia sets of transcendental entire functions has 
   also been studied by a number of authors (see e.g.~\cite{morosawa,bergweilermorosawa,osborne}). 
   This problem is closely connected to the boundedness of Fatou components, mentioned above. 
    Indeed, 
    if any component of $F(f)$ is unbounded,  then $J(f)$ cannot be locally connected 
   (compare \cite[Theorem E]{bakdom2000} and Lemma \ref{lem:new} below); 
   in particular, hyperbolic functions with asymptotic values  do not have locally connected Julia sets. 
    We may ask whether the conditions in Theorem \ref{thm:alljordan}, which describe precisely when 
     all Fatou components are bounded, also ensure local connectivity of the Julia set. It turns out that this is not the case:
\begin{example}[Non-locally connected Julia set]\label{ex:largedegree}
\label{example2}
There exists an entire function $f$, having two critical values~$0$ and~$1$
  and no other singularities of the inverse function, such that the following hold:
     \begin{enumerate}[{\normalfont(a)}]
   \item $0$ and $1$ are superattracting fixed points;
   \item every Fatou component of $f$ is bounded by a Jordan curve;
   \item the Julia set of $f$ is not locally connected.
     \end{enumerate}
\end{example}

The basic idea behind the construction is to use critical points of extremely high degree to simulate the behaviour of
  an asymptotic value and create (pre-periodic) Fatou components of large diameter; compare Figure \ref{fig:maclanevinberg-5-25}
   in Section \ref{sec:counterexamples}. This
  suggests that, to ensure local connectivity,  we should require a bound on the multiplicities of critical points within any one Fatou component. Indeed,
  using a result of Morosawa (see Theorem \ref{thm:morosawalc} below), we obtain the following consequence of
   Theorem \ref{thm:alljordan}.

\begin{cor}[Bounded degree implies local connectivity] \label{cor:lc}
  Let $f\in\B$ be hyperbolic with no asymptotic values. Suppose that there is a number $N$ 
   such that
   every component of $F(f)$ contains at most $N$ critical points, counting multiplicity. Then
   $J(f)$ is locally connected.
\end{cor}

Again, the additional assumption becomes particularly simple when every Fatou component contains at most one critical value, or when $\# \sing(f^{-1})=2$.
\begin{cor}[Locally connected Julia sets]\label{cor:lctwo}
Let $f$ be a hyperbolic function without asymptotic values. Suppose that 
\begin{enumerate}[\normalfont(a)] 
\item \label{item:lctwo_onecrit}every component of $F(f)$ contains at most one critical value, or
\item \label{item:lctwo}$\# \sing(f^{-1})=2$ and every component of $F(f)$ is bounded.
\end{enumerate}
Assume additionally that  the multiplicity of the critical points of $f$ is uniformly bounded. Then $J(f)$ is locally connected.
\end{cor}

 One interesting consequence of the preceding discussion is that, in the transcendental setting, \emph{local connectivity
    does not imply simple topological dynamics}. Indeed, the examples of Bishop that were mentioned above have  
   only critical points of degree at most $4$
    and, by precomposing such an example with a linear map, we can ensure that both critical values are superattracting fixed points. Hence 
    Corollary~\ref{cor:lctwo} applies to these families and therefore the Julia set is locally connected. On the other hand,
    by the results of \cite{rigidity}, the ``pathological'' behaviour near infinity is preserved by such a composition. In particular,
    we can find hyperbolic entire functions with locally connected Julia sets where the escaping set does not contain any curves to
    $\infty$, or even (using recent results from \cite{arclike}) such examples where the Julia set contains an uncountable
    collection of pairwise disjoint and dynamically natural pseudo-arcs. (The \emph{pseudo-arc} 
    is a certain hereditarily indecomposable continuum; cf.\ \cite[Exercise 1.23]{nadler}.)

\subsection*{Boundedness of immediate basins}
 The key step in establishing Theorem \ref{thm:alljordan} is to verify that all
    \emph{periodic} Fatou components of the map $f$ are bounded, provided that condition~\ref{item:allcritbounded}~in the theorem holds. 
    This is achieved by the following result, which gives a variety of conditions equivalent to the boundedness of a Fatou component. 
 \begin{thm}[Immediate basins of hyperbolic maps]\label{thm:jordan}
     Let $f\in\B$ be a hyperbolic transcendental entire function, and let $D$ be a periodic Fatou component of $f$, 
    say of period $p\geq 1$. 
   Then the following are equivalent:
 \begin{enumerate}[{\normalfont(a)}]
    \item $D$ is a quasidisc; \label{item:quasidisc}
    \item $D$ is a Jordan domain;\label{item:jordangeneral}
    \item $\Ch\setminus D$ is locally connected at some finite point of $\partial D$; \label{item:bakerweinreich}
    \item $D$ is bounded; \label{item:bounded}
      \item $D$ does not contain a curve to infinity;\label{item:curve}
     \item the orbit of $D$ contains no asymptotic curves and only finitely many critical points;
                  \label{item:critical}
     \item $f^p\colon D\to D$ is a proper map;\label{item:proper}
     \item for at least two distinct choices of $z\in D$, the set $f^{-p}(z)\cap D$ is finite. \label{item:finite}
 \end{enumerate}
  \end{thm}

As mentioned, the key new implication here is \ref{item:critical}$\Rightarrow$\ref{item:bounded}; the remaining equivalences and implications can be obtained
   by well-established methods, although some of them appear to be folklore. This part of the proof relies on a result from
   \cite{rigidity}, which states that hyperbolic maps are uniformly expanding on a suitable neighbourhood of their Julia sets; 
     see Proposition \ref{prop:expansion} below, and
    compare also
     \cite[Theorem C]{ripponstallardhyperbolic}. 

  We remark that the conclusion of the theorem does not hold if we omit the requirement that $f\in\B$ from the definition of
   hyperbolicity. Indeed, consider 
  $f(z) := e^z+z+1$, which is precisely Newton's method for   
  finding roots of $z\mapsto e^{-z}+1$. Then the singular values of $f$ 
  are precisely the infinitely many superattracting cycles $a_k=(2k+1)\pi i$ (with $k\in\Z$), and 
  $f$ has degree two when restricted to the invariant Fatou component containing $a_k$. However,
  all these components are unbounded.

 As a consequence of Theorem~\ref{thm:jordan}, we obtain the result announced after Theorem~\ref{thm:alljordan}, concerning quasidiscs:
\begin{cor}[Bounded components are quasidiscs]\label{cor:quasidiscs}
 Every bounded Fatou component of a hyperbolic entire function is a quasidisc.
\end{cor}
\begin{remark}[Remark 1]
 It is possible for pre-periodic unbounded Fatou components to be quasidiscs; indeed this is the case for two of the components 
   in Figure~\ref{fig:asymptotic}\subref{subfig:maclanevinberglimit}.
\end{remark}
\begin{remark}[Remark 2]
  This Corollary can also be deduced directly from known results and methods. Indeed, a theorem of Morosawa \cite[Theorem 1]{morosawa} 
    implies that every bounded Fatou component of a hyperbolic function is a Jordan domain. 
    It is not difficult to deduce that the boundary must in fact be a quasicircle.
\end{remark}

\subsection*{Bounded Fatou components and local connectivity beyond the hyperbolic setting}
  In this article, we consider only hyperbolic functions, and use uniform expansion estimates
  to establish our results. There are a number of weaker 
  hypotheses that will also suffice; here we mention only that all our 
   proofs go through for
  entire functions without asymptotic values that are \emph{strongly subhyperbolic} in the
  sense of Mihaljevi\'c-Brandt \cite{helenaorbifolds}. The theorems on Jordan Fatou components
  should extend even more generally, e.g.\ assuming that the function is \emph{geometrically finite}
   in the sense of \cite{helenalanding}, and there are no asymptotic values on the boundaries of
    Fatou components.
   (However, in the presence of parabolic points the boundaries will no longer be quasicircles.)
   In view of
    the recent result of Roesch and Yin \cite{roeschyin}  that all bounded
    attracting (and parabolic) Fatou components of polynomials are Jordan domains,
    it is plausible that the same always
    holds also in the entire transcendental setting, without additional 
    dynamical assumptions on the function $f$:

   \begin{conj}
      Let $f$ be a transcendental entire function, and let $D$ be an immediate attracting or
       parabolic basin. If $D$ is bounded, then $D$ is a Jordan domain. 
   \end{conj}

\subsection*{Structure of the article} In Section \ref{sec:preliminaries}, we shall collect 
  the prerequisites required to prove our theorems. The proof of Theorem
   \ref{thm:jordan} relies crucially on a uniform expansion estimate (Proposition 
   \ref{prop:expansion})
    for hyperbolic 
    entire functions, but we shall require a number of additional results to 
   deduce our theorems as stated. To keep the article self-contained, 
  to emphasise the elementary nature of our arguments, and to provide 
   a convenient reference for future studies of hyperbolic functions, we provide proofs or
   sketches of proofs where appropriate. We do use results of Heins \cite{heins} and
   Baker-Weinreich \cite{bakerweinreich} without further comments on their proofs. However, we 
    emphasise that these are required only to state properties~\ref{item:bakerweinreich} 
   and~\ref{item:finite} of Theorem~\ref{thm:jordan} in as weak a form 
   as possible, rather than being essential to the remainder of the proof.

 In Section \ref{sec:jordan}, we prove our main result, Theorem \ref{thm:jordan}, and 
   deduce the remaining theorems stated in the introduction in Section \ref{sec:dottingtheis}. 
   Finally, Section \ref{sec:counterexamples} is dedicated to the construction of 
   Examples \ref{ex:cossqrt}, \ref{example1} and \ref{example2}, using a method of
   MacLane and Vinberg.

\subsection*{Acknowledgements} We would like to thank Mario Bonk, S\'ebastien Godillon and J\"orn Peter for
  interesting discussions about this work. We are also grateful to 
   Alexandre Dezotti, Toby Hall, Mary Rees and Yao Xiao for their useful comments on the manuscript, and to the
   referee for their careful reading and detailed comments.%

\section{Preliminaries}\label{sec:preliminaries}

 \subsection*{Notation}
  As usual, we denote the complex plane by $\C$, and the Riemann sphere by $\Ch$. 
    Throughout the article, $f$ will denote a transcendental entire function, usually
    belonging to the Eremenko-Lyubich class $\B$ as defined in~\eqref{eqn:classB}.
   Recall from the introduction that $S(f) := \overline{\sing(f^{-1})}$ denotes the 
   set of \emph{singular values} of $f$.

   If $A,B\subset\C$, the notation $A\Subset B$ (``$A$ is compactly contained in $B$'') will
    mean that $A$ is bounded and $\overline{A}\subset B$.   The interior of a set $A\subset \C$ 
    is denoted by $\interior(A)$.

  We refer to \cite{waltersurvey} for background on transcendental iteration theory, and to
     \cite{beardonminda} for background on hyperbolic geometry.

 \subsection*{Hyperbolicity and uniform expansion}
 
Hyperbolicity is a key assumption in our results. We recall here some important properties. While these are well-known, we are not
  aware of a suitable reference, and hence provide a detailed proof here for the reader's convenience. 

\begin{prop}[Hyperbolic functions]\label{prop:hyproperties}
Let $f:\C\to\C$ be a transcendental entire function. 
Then $f$ is hyperbolic if and only if $\P(f)\Subset F(f)$.
 
If $f$ is hyperbolic, then 
   $F(f)$ is a finite union of attracting basins, and every
  connected component of $F(f)$ is simply-connected.
   Furthermore, there is a compact set $K\subset F(f)$ such that $f(K)\subset \interior(K)$  and
     $S(f)\subset \interior(K)$. The set $K$ can be chosen as the closure of a finite union of
    pairwise disjoint Jordan domains with analytic boundaries, no two of which belong to the
    same Fatou component.
\end{prop}
\begin{proof} 
First suppose that $f$ is hyperbolic. To see that $f$ has 
only finitely many attracting basins note that $S(f)$ is compact, and that the Fatou components of $f$ form an open covering of $S(f)$ by assumption.  Hence only finitely many 
 Fatou components intersect $S(f)$. On the other hand, every periodic cycle of attracting
  Fatou components contains at least one singular value by Fatou's Theorem \cite[Theorem 7]{waltersurvey}, so we see  that the number of such cycles is finite. 
It follows, in particular, that the postsingular set $\P(f)$ is compact and 
  contained in the Fatou set.

 Every connected component of an attracting basin is simply-connected by the maximum
  principle. We next construct the set $K$. We just proved that there are only finitely many
   Fatou components that intersect $\P(f)$. To each such component $U$,
   we can associate
  an analytic Jordan domain $D(U)\Subset U$ with $\P(f)\cap U\subset D(U)$ and
    $f(D(U))\Subset D(\widetilde{f(U)})$, where $\widetilde{f(U)}$ denotes the Fatou component containing
    $f(U)$. (Compare \cite[Proposition~2.6]{helenaorbifolds}.)
    Indeed, we first construct the domains $D(U)$ for all pre-periodic components by induction, beginning with the
    components having the largest pre-period. Finally, let $\rho\gg 0$ be a sufficiently large constant, and for 
    each periodic Fatou component, let $D(U)$ be the hyperbolic disc of radius $\rho$ around the attracting periodic point in $U$.
    Then $f(D(U))\Subset D(\widetilde{f(U)})$ by the Schwarz lemma. If $\rho$ was chosen sufficiently large, then $D(U)$ also contains
    the compact set $\P(f)\cap U$ and the image of $D(V)$ for any pre-periodic component $V$ with $\P(f)\cap V\neq\emptyset$ and $f(V)\subset U$. 
  Since, by construction, $f(D(U))\Subset D(\widetilde{f(U)})$ for all components $U$ of pre-period greater than
  $1$ that intersect $P(f)$, we see that
   \[ K := \bigcup_{U\cap \P(f)\neq\emptyset} \overline{D(U)} \] 
   has the required properties.

From now on, assume only that $\P(f)\Subset F(f)$, 
   and let $U$ be a component of $F(f)$. Then $U$ cannot be a Siegel disc,
  since otherwise we would have $\partial U \subset \P(f)$ \cite[Theorem~7]{waltersurvey}, and hence $\P(f)\cap J(f)\neq \emptyset$, which contradicts
   our assumption.
This implies that all limit functions of the family $\left\{f^n|_U\colon n\in\N\right\}$
are constant, possibly infinite.
Since $f\in\B$, we cannot have $f^n|_U\to\infty$, 
by a result of Eremenko and Lyubich \cite[Theorem~1]{eremenko_lyubich_2}.
Hence there exists a subsequence of $\left(f^n|_U\right)$
which tends to a finite constant $a\in\C$.
A~result of Baker \cite[Theorem~2]{bak70} implies that $a\in \P(f)$.
By assumption, it follows that $a\in F(f)$. This implies that $U$ can be neither a parabolic nor a wandering domain.
It follows that $a$ must in fact be an attracting periodic point, and $U$ is a component of its attracting basin.
Hence $F(f)$ consists of only of attracting basins. In particular, $f$ is hyperbolic, and the proof is complete.
\end{proof}
We remark that Bishop \cite{bishopfolding} recently proved that    the class $\B$ contains 
  non-hyper\-bolic functions that do have wandering domains. (The orbits of these domains
   accumulate   both at $\infty$ and at some finite points.
Wandering domains with the latter property had been constructed
earlier by Eremenko and Lyubich~\cite[Example~1]{elexamples}, but in their examples it was 
not clear whether the function could be taken to be in~$\B$.)

The key element in our proof of Theorem \ref{thm:jordan} will be the fact that hyperbolic functions are
 uniformly expanding, with respect to a suitable conformal metric.
\begin{prop}[Uniform expansion for hyperbolic functions {\cite[Lemma 5.1]{rigidity}}]\label{prop:expansion}
  Let $f\colon\C\to\C$ be a hyperbolic transcendental entire function, and let $K$ be the compact set from 
   Proposition \ref{prop:hyproperties}. That is, $f(K)\subset\interior(K)$ and $S(f)\subset K$.

  Define $W := \C\setminus K$ and $V:= f^{-1}(W)$. Then there is a constant $\lambda>1$ such that
     \[ \|Df(z)\|_W \geq \lambda \]
    for all $z\in V$, where $\|Df\|_W$ denotes the derivative of $f$ with respect to the hyperbolic metric of $W$.
\end{prop}
\begin{proof}[Idea of the proof]
For completeness, let us briefly sketch the proof of this fact; we refer to \cite{rigidity} for details. Since $f\colon V\to W$
   is a covering map and $\cl{V}\subset W$, we have 
    \[   \|Df(z)\|_W = \frac{\rho_V(z)}{\rho_W(z)}>1 \]
  for all $z\in V$. It follows that it suffices to prove
     \[ \rho_W(z) = o(\rho_V(z)) \]
   as $z\to\infty$. By standard estimates, we have
      \[ \rho_W(z) = O\left(\frac{1}{|z|\log|z|}\right), \]
    while for $V$ it is shown in \cite{rigidity} that
       \[ \frac{1}{\rho_V(z)} = O(|z|). \]
     (This uses the fact that $\C\setminus V = f^{-1}(K)$ contains a sequence $(w_n)$ with    
         $|w_{n+1}| \leq C|w_n|$, for a constant $C>1$, together with estimates on the hyperbolic metric in a
    multiply-connected domain.) This completes the proof.
\end{proof}

\subsection*{Local connectivity}

We shall use the following characterization of local connectivity for compact subsets of the Riemann sphere.
\begin{lem}[LC Criterion~{\cite[Thm. 4.4, Chapter VI]{whyburn1942}}]\label{whyburnlemma}\label{lem:whyburn}
A compact subset of the Riemann sphere is locally connected if
and only if the following two conditions are satisfied:
\begin{enumerate}[{\normalfont(a)}]
\item the boundary of each complementary component is locally connected;\label{item:whyburnlc}
\item for every positive $\varepsilon$ there are only finitely many complementary 
components of spherical diameter greater than~$\varepsilon$.\label{item:whyburndiameter}
\end{enumerate}
\end{lem}
We will apply the result above to study the local connectivity of Julia sets, 
even though we consider these to be subsets of $\C$ rather than of $\Ch$. However, it is well-known that a continuum cannot fail to
  be locally connected at a single point \cite[Corollary~5.13]{nadler}. Hence it
  follows that $J(f)$ is locally connected if and only if $J(f)\cup\{\infty\}$ is.

As mentioned in the introduction, a locally connected Julia set cannot have an unbounded Fatou component
 \cite[Corollary~1.2]{osborne}.
\begin{lem}[Local connectivity implies bounded Fatou components] 
\label{lem:new}
Let $f$ be an entire transcendental function. If $F(f)$ has an unbounded Fatou component, then $J(f)$ is not locally connected.
\end{lem}
\begin{proof}
Let $U$ be an unbounded component of $F(f)$. First suppose that there is some iterated preimage component $\tilde{U}$ of $U$
  that is not periodic. It follows that, for any bounded open set $D$ intersecting the Julia set, 
 there are infinitely many different unbounded Fatou components (namely iterated preimages of $\tilde{U}$) 
 that intersect $D$. Hence condition \ref{item:whyburndiameter} of Lemma \ref{whyburnlemma} is violated, and $J(f)$ is not locally connected.

If no such component $\tilde{U}$ exists, then $U$ is completely invariant for some iterate $f^n$. It follows that 
   $J(f)=J(f^n)$ is not locally connected by  \cite[Corollary 3]{bakdom2000}. 
\end{proof}

 The following result shows that, once we know that all immediate attracting basins are  Jordan,   we can already make conclusions about the local connectivity of the Julia set~-- provided   that there is a bound on the degree of $f$ on any pre-periodic Fatou component. 

 \begin{thm}[Bounded components and bounded degree imply local connectivity]\label{thm:morosawalc}
 Let $f\in \B$ be hyperbolic with no asymptotic values. Suppose that
   every immediate attracting basin of $f$ is a Jordan domain. If there is $N$ such
    that the degree of the restriction of $f$ to any Fatou component is bounded by $N$, then
    $J(f)$ is locally connected.
 \end{thm}
   This theorem is due to Morosawa \cite[Theorem 2]{morosawa}. We note, however, that the
     statement in \cite{morosawa} overlooked the assumption on the degree of
     preimages of Fatou components. Our Examples \ref{ex:unbounded} and \ref{ex:largedegree}
     show that this assumption is necessary. A more general statement (which includes the corrected hypotheses) can be found in \cite[Theorem 4]{bergweilermorosawa}. For convenience, 
    let us show how the
     result can be obtained from Proposition \ref{prop:expansion}.
 \begin{proof}[Proof of Theorem \ref{thm:morosawalc}]
    We only need to establish part \ref{item:whyburndiameter} of Lemma \ref{lem:whyburn}. 
     By Proposition~\ref{prop:hyproperties}, all components of $F(f)$ are simply-connected.
     Hence, if $U$ and $V$ are Fatou components with $f(U)\subset V$ 
     and $V\cap S(f)=\emptyset$, then $f\colon U\to V$ is bijective.
     Since $S(f)$ is compactly contained in the Fatou set,
    only a finite number $k$ of Fatou components intersect the singular set.
     
Let $U$ be a pre-periodic Fatou component, say of
       pre-period $n$, and let $V=f^n(U)$ be the first periodic Fatou component on the orbit of $U$.
  The assumption  implies that the degree of $f^n\colon U\to V$ is bounded by $N^k$.
   Hence ~-- using that $V$ is a Jordan domain~-- 
    the boundary $\partial U$ covers $\partial V$ at most $N^k$
    times when mapped under $f^n$. 
    Let $W$ and $\lambda$  be as in Proposition \ref{prop:expansion}. We can cover $\partial V$ by, say, $M$ simply-connected hyperbolic discs 
     (with respect to the hyperbolic metric of $W$). Since $f$ has only finitely many
    periodic Fatou components, the number $M$ is bounded independently of $U$.
    Let $r$ be the maximal hyperbolic diameter
of these discs. Proposition   \ref{prop:expansion} implies that we can cover $\partial U$ by $N^k M$  hyperbolic disks of diameter $r/\lambda^n$. 
Thus the hyperbolic diameter of $U$ in $W$ is bounded by
$N^k Mr/\lambda^n$,  which tends to zero exponentially as $n$ tends to infinity. Furthermore, for a given $n$, the Fatou components of pre-period $n$ can only accumulate at infinity (by the open mapping theorem), and hence only finitely
many of them have spherical diameter greater than a given number $\eps>0$. This establishes property \ref{item:whyburndiameter} of Lemma
      \ref{lem:whyburn}, and completes the proof.
\end{proof}
\begin{remark}
   Alternatively, we could use distortion principles for maps of bounded degree to see that every Fatou component $U$ contains
    an open disc of size comparable to the diameter of $U$. 
  Again, this implies~\ref{item:whyburndiameter} in Lemma \ref{lem:whyburn}.
 \end{remark}

In order to state some of our results concerning local connectivity of Fatou component boundaries, we shall use the following theorem,   
   which is due to Baker and  Weinreich \cite{bakerweinreich}. We remark that this result is not central to our arguments, but rather allows
    us to state some conclusions (e.g.\ in Theorem \ref{thm:jordantwo}) more strongly than would otherwise be possible.
\begin{thm}[Boundaries of periodic Fatou components]\label{thm:bakerweinreich}
 Let $f$ be a transcendental entire function, and suppose that $U$ is an unbounded periodic component of
    $F(f)$ such that $f^n|_U$ does not tend to infinity. Then $\Ch\setminus U$ is not locally connected at any finite point of
       $\partial U$.
\end{thm}
\begin{proof}
 Baker and Weinreich proved that, under these assumptions, the impression of every prime end of $U$ contains $\infty$. 
   Equivalently, if $\phi\colon \D\to U$ is a conformal map (which exists by the Riemann mapping theorem) and $z_0\in\partial \D$, then there exists a sequence $z_n\in\D$ such that $z_n\to z_0$ and $\phi(z_n)\to \infty$.

Now suppose, by contradiction, that some point in $\partial U$ has a bounded connected neighbourhood $K$ in $\Ch\setminus U$, which we may assume to
   be compact and full. Let $z_0\in\partial \D$ such that the radial limit of $\phi$ at
    $z_0$ exists and belongs to the relative interior of $K$ in $\Ch\setminus U$. 
 There is a small round disc $D$ around $z_0$ such that 
    the Euclidean length of $\gamma := \phi(\D\cap \partial D)$ is sufficiently short to ensure that   
   both endpoints of $\gamma$ are in $K$, and 
    that $K\cup\gamma$ does not separate $\phi(0)$ from $\infty$. (This follows from 
    a well-known application of
     the length-area principle~-- see e.g.\ \cite[Proposition~2.2]{pommerenkeboundary}, which strictly 
     speaking 
     applies only to bounded domains, but whose proof 
    yields the desired result in the unbounded case
     upon replacing Euclidean length and area
    with their spherical analogues.) 
  It follows that $\phi(D\cap\partial \D)$ is contained in the bounded complementary 
    component of $K\cup\gamma$, which is a contradiction to the above result by Baker and Weinreich.
    (Alternatively, we may additionally assume that the prime end corresponding to $z_0$ is symmetric, as the set of
     asymmetric prime ends is at most countable \cite[Proposition 2.21]{pommerenkeboundary}. By 
     \cite[Corollary 1]{carmonapommerenke}, the corresponding impression  is contained in $K$~-- a contradiction.)
\end{proof}

Finally,  we shall require a number of facts concerning the mapping behaviour
  of entire functions on preimages of simply-connected domains. While these results are certainly  not new, we are again not aware of a convenient reference and therefore include the proofs.
In our arguments, we shall use the following simple lemma.
\begin{lem}[Coverings of doubly-connected domains]\label{lem:annulus-covering}
Let $A,B\subset\C$ be domains and let $f\colon B\to A$ be a covering map. 
Suppose that $A$ is doubly-connected. Then either $B$ is doubly-connected and $f$ is a 
proper mapping, or $B$ is simply-connected (and $f$ is a universal cover, of infinite degree).
\end{lem}
\begin{proof}
The fundamental group of $A$ is isomorphic to $\Z$.
The fundamental group of $B$ is thus isomorphic to a subgroup of $\Z$.
As the only subgroups of $\Z$ are the trivial one and the groups $k\Z$ with
$k\geq 1$, the conclusion follows easily.
\end{proof}

\begin{prop}[Mapping of simply-connected sets]\label{prop:inner}
 Let $f$ be an entire function, let $D\subset\C$ be a simply-connected domain, and
   let $\widetilde{D}$ be a component of $f^{-1}(D)$. Then either
\begin{enumerate}[{\normalfont(a)}]
  \item \label{item:inner_proper}
    $f\colon\widetilde{D}\to D$ is
   a proper map (and hence has finite degree), or 
  \item \label{item:inner_infinite}
    $f^{-1}(w)\cap \widetilde{D}$ is infinite
     for every $w\in D$, with at most one exception. 
\end{enumerate}
 In case~\ref{item:inner_infinite}, either $\widetilde{D}$ contains an asymptotic curve corresponding
  to an asymptotic value in $D$, or $\widetilde{D}$ contains infinitely many critical points.
\end{prop}
\begin{proof}
 A theorem by Heins \cite[Theorem 4']{heins} implies that either~\ref{item:inner_infinite} holds, or the number of preimages of $w\in D$ in $\widetilde{D}$
  (counting multiplicity) is finite and constant in $D$. It is elementary 
  to see that the latter is equivalent 
  to~\ref{item:inner_proper}.

 To prove the final statement, it is sufficient to consider the case where 
   $f\colon\widetilde{D}\to D$ has no asymptotic values in $D$ and only finitely many critical values (otherwise, there is nothing to show).
  This implies that this map is an infinite \emph{branched covering}; i.e.,  
   every point $z_0\in D$ has a simply-connected neighbourhood $U$ such that every 
   component $\widetilde{U}$ of $f^{-1}(U)\cap D$ is mapped as a finite covering, branched
   at most over $z_0$. 
 
 Such a map must have infinitely many critical points. This essentially follows from the
  Riemann-Hurwitz formula~-- which is usually stated only for proper maps, but whose
  proof goes through also in this case. For completeness, let us 
   indicate an alternative proof of our claim. 
Let $c_1\dots,c_m$ be the distinct critical values in $D$ 
%
We join them to a further point $a\in D$ by simple arcs $\tau_1,\dots,\tau_m$
which do not intersect except in their common endpoint~$a$.
Set $T=\bigcup_{j=1}^m \tau_j$ 
and $\widetilde{T}:=f^{-1}(T) \cap \widetilde{D}$. 
  Now $D\setminus T$ is doubly-connected and $f\colon \widetilde{D} \setminus \widetilde{T} \to  D\setminus T$ is a covering of infinite degree. 
By Lemma~\ref{lem:annulus-covering} this map 
    must be a universal covering,  and thus 
     every component $T'$ of $\widetilde{T}$ is unbounded. 
     Since $f$ is a branched covering map, $T'$ consists of infinitely many preimages of $T$, joined
     together only at critical points. This 
     proves that $\widetilde{D}$ contains infinitely many critical points.
 \end{proof}

In our applications, the function $f$ will always be hyperbolic, and hence the
  set of singular values stays away from the boundary of $D$. In this case, we can say
  more:

\begin{prop}[Preimages of sets with non-singular boundary] \label{prop:covering}
  Let $f$, $D$ and $\widetilde{D}$ be as in Proposition \ref{prop:inner}, and assume
   additionally that $D\cap S(f)$ is compact.
  \begin{enumerate}[{\normalfont(1)}]
   \item If $\# D\cap S(f)\leq 1$, then $\widetilde{D}$ contains at most one critical point 
  of $f$. \label{item:onecriticalpoint}
   \item If $S(f)\subset D$, then $\widetilde{D}=f^{-1}(D)$. \label{item:connectedpreimage}
   \item In case \ref{item:inner_proper} of Proposition \ref{prop:inner},
     if $D$ is a bounded Jordan domain (resp.~quasidisc) such that
     $\partial D\cap S(f)=\emptyset$, then $\widetilde{D}$ is also a bounded Jordan domain (resp.\ quasidisc).\label{item:jordaninduction}
   \item In case \ref{item:inner_infinite} of Proposition \ref{prop:inner}, the point $\infty$ is accessible from $\widetilde{D}$.\label{item:criticalaccessible}
  \end{enumerate}
\end{prop}
\begin{proof}
Set $S_D := S(f)\cap D$. If $\# S_D = 1$, then 
    $f\colon\widetilde{D}\setminus f^{-1}(S_D) \to D\setminus S_D$ is conformally equivalent to an unramified covering of
    the punctured unit disc. By Lemma \ref{lem:annulus-covering}, it follows that $\widetilde{D}$ contains at most one critical point.
We have proved~\ref{item:onecriticalpoint}.

Now let $U\subset D$ be simply-connected such that  
     $S_D\subset U \Subset  D$. Set $\widetilde{U} := f^{-1}(U)\cap \widetilde{D}$. 
	     By the maximum principle, every component of $\widetilde{U}$ is simply-connected. We will show  that 
     $\widetilde{U}$ is connected. 
   Indeed, let $z,w\in \widetilde{U}$, and let $\widetilde{\gamma}\subset \widetilde{D}$ 
    be a smooth arc connecting $z$ and $w$.  
    Set $\gamma=f(\widetilde{\gamma})$. Since $S_D$ has a positive distance from $\partial U$,
the curve $\gamma$ contains at most finitely pieces that connect $S_D$ to $\partial U$. By
    cutting the curve at one point in each of these pieces, we may divide it into
    finitely many segments, some that may intersect $S(f)$ but are contained in $U$, and some
    that may leave $U$  but do not intersect $S(f)$. We construct a new curve   
    $\gamma' \in U$ which equals $\gamma$ in those segments contained in $U$ 
    but is only homotopic to $\gamma$ in $D\setminus S(f)$, relative to the endpoints, in the 
   remaining pieces. These homotopies can be lifted to $\widetilde{D}\setminus f^{-1}(S(f))$ 
   since $f$ is a covering there, resulting in a curve $\widetilde{\gamma}' \subset \widetilde{U}$ 
   connecting $z$ and $w$. Hence $\widetilde{U}$ is connected.
   In particular, this proves~\ref{item:connectedpreimage} 
   (replacing $D$ by $\C$ and $U$ by $D$).

  For the remainder of the proof, let us require additionally that
     $U$ is a bounded Jordan domain with $U\Subset D$. 
     Then $A:= D\setminus U$ is doubly-connected. Consider the set $\widetilde{A} := f^{-1}(A)\cap \widetilde{D}$. On every component of $\widetilde{A}$, the restriction of 
     $f$ is a holomorphic covering map, since $S_D\subset U$. 
By Lemma~\ref{lem:annulus-covering}, the components of $\widetilde{A}$ are
either doubly- or simply-connected.
      
  Suppose first that $\widetilde{A}$ has a doubly-connected component. Since $\widetilde{U}$ is connected 
    and simply-connected, it follows that
      $\widetilde{A}$ is connected, and that $\widetilde{U}$ is bounded. As 
      $f\colon \widetilde{A}\to A$ has finite degree, it follows that we are in case
      \ref{item:inner_proper} of Proposition \ref{prop:inner}. 
       Hence $f\colon\widetilde{D}\to D$ is a proper map. 
      If, additionally, $D$ is a bounded Jordan domain 
     whose boundary is disjoint from $S(f)$, then 
       we can apply Proposition \ref{prop:inner} to a slightly larger Jordan disc 
       $D'$ without additional singular values, so 
the restriction of $f$ to the preimage of $D'$ containing $\widetilde{D}$ is still proper.
Then, $f\colon \partial \widetilde{D} \to \partial D$ is a finite degree covering map,
      which proves that $\partial \widetilde{D}$ is indeed a Jordan domain. 
    Of course the property of being a quasicircle is preserved under a conformal covering map.  
    This establishes \ref{item:jordaninduction}.
       
Now suppose that every component $V$ of $\widetilde{A}$ is simply-connected. 
   Then $f\colon V\to A$ is a universal covering, and hence has infinite degree. 
     The preimage of any simple 
    non-contractible closed curve in $A$ under this covering is a Jordan arc in $\widetilde{A}$ tending to
     infinity in both directions, and hence $\infty$ is accessible from $\widetilde{D}$,
   proving \ref{item:criticalaccessible}.
\end{proof}
\begin{remark}[Remark 1]
 In \ref{item:jordaninduction}, to conclude that $\tilde{D}$ is bounded,
   it is enough to assume that $\partial D$ has exactly two complementary components,
   rather than that $\partial D$ is a Jordan curve. Indeed, this follows from the original statement,
   since we can surround $\overline{D}$ by a Jordan curve $\gamma$
   such that the Jordan domain $W$ bounded by $\gamma$ does not contain any 
   singular values other than those already in $D$. The claim then follows from the Proposition
   as stated.
\end{remark}

\begin{remark}[Remark 2]
In the case when $f\colon\widetilde{D} \to D$ is proper but $\widetilde{D}$ is unbounded, we do not know
   whether $\infty$ is always accessible from $\widetilde{D}$. Indeed, this is an open question even
   when $D=\widetilde{D}$ is an unbounded Siegel disc of an exponential map. Also compare
   the question in \cite[p.\ 439, ll.\ 8--9]{BakDom99}.
\end{remark}

\section{Periodic Fatou components}\label{sec:jordan}

\begin{proof}[Proof of Theorem \ref{thm:jordan}]
 Let $f\in\B$ be hyperbolic, and let $D$ be an immediate attracting basin of $f$, say of period $p$.
 By passing to an iterate, we may assume without loss of generality that $p=1$.
  Recall that $D$ is simply-connected (by the maximum principle).

  Clearly \ref{item:quasidisc}$\Rightarrow$\ref{item:jordangeneral}$\Rightarrow$\ref{item:bakerweinreich}, as
    any quasidisc is Jordan, and the complement of any Jordan domain is locally connected at every point.   On the other hand, if $\Ch\setminus D$ is locally connected at \emph{any} finite point of $\partial D$, then $D$ is bounded by
    Theorem \ref{thm:bakerweinreich}. So \ref{item:bakerweinreich} implies 
     \ref{item:bounded}.

   Clearly, if $D$ is bounded, then $D$ cannot contain a curve to $\infty$, and hence
     \ref{item:bounded}$\Rightarrow$\ref{item:curve}. 

  Since $f$ is hyperbolic, $S(f)\cap D$ is compact and we may apply 
   Proposition~\ref{prop:covering} \ref{item:criticalaccessible} to conclude that,
    if $D$ does not contain a curve to infinity, then alternative~\ref{item:inner_proper} of Proposition~\ref{prop:inner}
holds. This in turn implies that $D$ contains only finitely many critical points and no asymptotic values, and,  again by 
     Proposition \ref{prop:inner}, this is in turn implies that $f\colon D\to D$ is a proper map. Hence 
    \ref{item:curve}$\Rightarrow$\ref{item:critical}$\Rightarrow$\ref{item:proper}.
     Since any proper map has finite degree, we have 
     \ref{item:proper}$\Rightarrow$\ref{item:finite}.

It remains to prove \ref{item:finite}$\Rightarrow$\ref{item:quasidisc}.
    So suppose that at least two points in $D$ each have
      at most finitely many preimages in $D$. We must show
     that $\partial D$ is a quasicircle.
    We shall first prove that
      $\partial D$ is a bounded curve. In the rational case, this argument goes back to 
     Fatou \cite[p.\ 83]{fatoumemoir2}; compare
     \cite[Chapter 5, Section 5]{steinmetz}, and using the uniform expansion from Proposition \ref{prop:expansion}, the proof goes
      through essentially verbatim. 

  To provide the details, let $U_0\Subset D$ be a bounded Jordan domain with analytic boundary such that
    $\overline{f(U_0)}\subset U_0$ and $S(f)\cap D\subset U_0$; such a domain exists by Proposition~\ref{prop:hyproperties}.
  By Proposition \ref{prop:inner}, $f\colon D\to D$ is 
    a proper map of some degree $d\geq 1$. 
   For $n\geq 1$, set 
    $U_n:= f^{-n}(U_0)\cap D$. Then $D=\bigcup {U_n}$.

  Now, due to the choice of $U_0$, we see that 
    $f^n\colon D\setminus \overline{U_n}\to A := D\setminus\overline{U_0}$
    is a finite-degree covering map (of degree $d^n$) over the doubly-connected domain $A$. 
    By Lemma~\ref{lem:annulus-covering}, the domain
    $D\setminus\overline{U_n}$ is also doubly-connected, and hence $U_n$ is connected for 
    all $n$.
    Furthermore, by Proposition~\ref{prop:covering}~\ref{item:jordaninduction}, applied to 
     $U_n$ and $U_0$, we see that
    each $U_n$ is a Jordan domain. Hence
    $f\colon \partial U_{n+1}\to \partial U_n$ is topologically a $d$-fold covering over a circle
     for every $n\geq 0$.

    We claim that we can find a diffeomorphism
    \[ \phi \colon W\to D\setminus \overline{U_0}, \qquad  W := \{z\in\C\colon 1/e<|z|<1\}, \]
   such that    
     \begin{equation}\label{eqn:functionalequationphi}
       f(\phi(z)) = \phi(z^d) \qquad \text{when } e^{-1/d}<|z|<1. \end{equation}

   Indeed, since $f\colon \partial U_1\to \partial U_0$ is a $d$-fold covering, we can
    define $\phi$ on
    the circles $\{\ln |z|=-1\}$ and $\{ \ln|z|= -1/d\}$ so that
    the functional relation~\eqref{eqn:functionalequationphi} is satisfied. By interpolation, we extend $\phi$ to a diffeomorphism
    \[ \left\{ z\in\C\colon -1 < \log |z|\ < -1/d\right\}\to A_0 := U_1\setminus\overline{U_0}. \] 

  Consider the annuli $A_n := U_{n+1}\setminus \overline{U_n}$. Since each $f\colon A_{n+1}\to A_n$ is a $d$-fold covering map of annuli, we can 
    inductively lift $\phi$ to a map 
       \[ \left\{ z\in\C\colon -d^{-n} < \ln |z| < -d^{-(n+1)} \right\} \to A_n, \]
    and our initial choice of $\phi$ ensures that this lift can be taken to extend the original map continuously. This completes the construction
    of $\phi$.

   Now let $\theta\in\R$ and $n\geq 0$, and consider the curve
      \[ \gamma_{n,\theta} := \phi( \{ e^{a +i\theta} \colon -d^{-n} \leq a \leq -d^{-(n+1)}\} ). \]
   By the functional relation~\eqref{eqn:functionalequationphi}, $\gamma_{n,\theta}$ is the image of the arc
   $\gamma_{0,\theta\cdot d^n}$ under some branch of $f^{-n}$. Recall from Proposition \ref{prop:expansion} that 
      \begin{equation}\label{eqn:expansionlambda} \|Df(z)\|_W \geq \lambda \end{equation}
    whenever $z,f(z)\notin K$, where $K$ is the compact set from Proposition~\ref{prop:hyproperties},
    $W=\C\setminus K$, and $\lambda>1$ is a suitable constant. 

  Since $D\setminus U_0\subset W$, it follows that
     \[ \ell_W(\gamma_{n,\theta}) \leq \lambda^{-n} \ell_W(\gamma_{0,\theta\cdot d^n}) \leq
          \lambda^{-n} \max_{\widetilde{\theta}} \ell_W(\gamma_{0,\widetilde{\theta}}). \]
 
  Thus, for $n\geq 0$, the functions 
  \[
  \sigma_n\colon \R/\Z \to D \qquad \sigma_n(t)=\varphi\left(e^{-d^{-n}+ 2\pi i t}\right) \in \partial U_n
  \]
  form a Cauchy sequence in the hyperbolic metric of $W$ as $n\to \infty$. Hence there exists a limit function $\sigma_n\to \sigma$ which is the continuous extension of $\varphi$ to the unit circle.   It follows that
     $\partial D$ is indeed a continuous closed curve. By the maximum principle,
     $\partial D = \partial \overline{D}$ and $\C\setminus \overline{D}$ has no bounded connected components. Hence 
     $\partial D$ is a Jordan curve.

  To see that $\partial D$ is a quasicircle, we again use the expanding property~\eqref{eqn:expansionlambda} of $f$ to find a
    Jordan neighbourhood $\Omega$ of $\overline{D}$ such that
    $f\colon\Omega\to f(\Omega)$ is a branched covering map of degree $d$ with 
  $\overline{\Omega}\subset f(\Omega)$.  
    We can now apply the Douady-Hubbard straightening theorem 
     \cite[Theorem 1, p.\ 296]{douhub85} 
    to see that $f|_{\Omega}$ is quasiconformally conjugate to a hyperbolic polynomial of degree 
   $d$ with an attracting fixed point whose immediate basin contains all the critical values. Such 
   a basin is completely invariant under the polynomial and its boundary is a quasicircle;
   hence $\partial D$ has the same properties.
   \end{proof}
   \begin{remark}[Remark]
  Alternatively, to  avoid 
    the use of the straightening theorem, it is possible to
    verify the geometric definition of quasicircles directly. This definition requires that the diameter
    of any arc of $\partial D$ is comparable to the distance between its endpoints. That
    condition is trivially satisfied on big scales, and we can transfer it to arbitrarily small
    scales using univalent iterates. (We thank Mario Bonk for this observation.)
 \end{remark}
 
\section{Bounded Fatou components and local 
   connectivity of Julia sets}\label{sec:dottingtheis}

 We now deduce the remaining theorems stated in the introduction, using Theorem~\ref{thm:jordan}.

 \begin{proof}[Proof of Corollary \ref{cor:quasidiscs}]
  Let $f$ be hyperbolic. 
   By Theorem \ref{thm:jordan}, any bounded \emph{periodic} component of $F(f)$ is a quasidisc. 
   Now, if $U$ is any bounded Fatou component, then clearly $U$ contains only finitely many critical points and no asymptotic 
   curves. Hence $f\colon U\to f(U)$ is a proper map by 
   Proposition \ref{prop:inner}, and if $f(U)$ is a quasidisc, then $U$ is also a quasidisc by  
   part ~\ref{item:jordaninduction} of
    Proposition \ref{prop:covering}. Hence, by induction every bounded Fatou component of $f$ is a quasidisc, as claimed.
 \end{proof}

 \begin{proof}[Proof of Theorem \ref{thm:alljordan}]
  Let $f\in\B$ be hyperbolic. If $f$ has no asymptotic values and every Fatou component
    contains at most finitely many critical points, then every periodic Fatou component is a bounded quasidisc 
    by Theorem~\ref{thm:jordan}.
    Moreover, by Proposition \ref{prop:inner} and part ~\ref{item:jordaninduction} of
    Proposition \ref{prop:covering},
      if $U$ is any Fatou component of $f$, then $f\colon U\to f(U)$ is a proper map, and if $f(U)$ is a bounded quasidisc, then so is $U$. 
      By induction on the pre-period of $U$, it follows all Fatou components are indeed bounded quasidiscs. 

   On the other hand, if $f$ has an asymptotic value, this value belongs
    to the Fatou set by hyperbolicity. Hence $f$ has an unbounded Fatou component. 
    Finally if  some Fatou component $U$ contains infinitely many critical points, these can only accumulate at infinity and therefore $U$ is unbounded.  
\end{proof}

\begin{proof}[Proof of Corollary \ref{cor:onlyonecriticalpoint}] 
  Let $f$ be hyperbolic with no asymptotic values, and assume that every Fatou component contains at most
    one critical value. By Proposition \ref{prop:covering}~\ref{item:onecriticalpoint}, it follows that 
     each Fatou component also contains at most one (possibly high-order) critical point. Thus every Fatou component is bounded
     by Theorem \ref{thm:alljordan}. 
\end{proof}

\begin{proof}[Proof of Theorem \ref{thm:jordantwo}]
  Suppose that $f$ is hyperbolic without asymptotic values, and with exactly two critical values. 
    Assume first that both critical values belong to the same Fatou component $D$.
    Then $D_0 := f^{-1}(D)$ is connected by Proposition \ref{prop:covering} \ref{item:connectedpreimage} and unbounded by the
     Casorati-Weierstrass theorem. Thus  $D_0$ is an unbounded component of $F(f)$ that contains
    all critical points of $f$. 
     By Fatou's theorem \cite[Theorem 7]{waltersurvey}, each cycle of attracting periodic components of $F(f)$ must contain a critical point, and hence 
     $D_0$ is periodic. 
    By Theorem \ref{thm:bakerweinreich},  $\partial D_0$ is not locally connected at any point. 
     
     Moreover, if $U$ is any component of $F(f)$, then, by Proposition
     \ref{prop:hyproperties}, there exists a minimal number
     $k$ such that $f^k(U)=D_0$. 
    Since $\infty$ is accessible from $D_0$ by  
     Proposition~\ref{prop:covering}~\ref{item:criticalaccessible}, this implies
that $\infty$ is accessible from every Fatou component, including $D$. In particular, all these components are unbounded.
 In order to prove that $\partial U$ is not locally connected at any finite point,
     we shall show that $f^k$ maps $\partial U$
     homeomorphically to $\partial D_0$. Indeed, let $\Gamma\subset D$ be an arc to infinity that  
    contains both critical values. 
     By the choice of $k$, we know that $f^j(U)$ is disjoint from $D$ for $j=1,\dots,k$, since $D$ has no preimage components
    apart from $D_0$.
    Thus there is a branch of $f^{-1}$ on
     $\C\setminus\Gamma$ that maps $\overline{f^k(U)}$ to $\overline{f^{k-1}(U)}$ homeomorphically, and hence we have established 
     case \ref{item:completelyinvariant} of Theorem~\ref{thm:jordantwo}.

   Assume now that the critical values are not both in the same
     Fatou component.
     Then  every Fatou component contains at most one  of them.
     By Corollaries \ref{cor:onlyonecriticalpoint} and \ref{cor:quasidiscs}, case~\ref{item:jordan} of Theorem
    \ref{thm:jordantwo} is satisfied. 
\end{proof}

 \begin{proof}[Proof of Corollary \ref{cor:lc}]
  Let $f\in\B$ be hyperbolic with no asymptotic values, let $N\in\N$, 
    and suppose that every Fatou component $U$
    of $f$ contains at most $N$ critical points (counting multiplicity). By Theorem \ref{thm:alljordan}, 
    every Fatou component $U$ is a bounded quasidisc,
    and the restriction $f\colon U\to f(U)$ is a proper map (see the proof of Corollary~\ref{cor:quasidiscs}). 
    By the Riemann-Hurwitz formula, the degree of this restriction is bounded 
     by $N+1$, since all components are simply-connected. Thus $J(f)$ is locally connected by Theorem~\ref{thm:morosawalc}.
 \end{proof}

\begin{proof}[Proof of Corollary \ref{cor:lctwo}]
  If~\ref{item:lctwo_onecrit} is satisfied, every Fatou component is bounded by Corollary \ref{cor:onlyonecriticalpoint}. 
  By hypothesis, the multiplicity of the critical points is uniformly bounded, and hence we may apply Corollary~\ref{cor:lc} and conclude that $J(f)$ is locally connected.
  In case~\ref{item:lctwo}, it was shown in the proof of Theorem \ref{thm:jordantwo}, that every Fatou component of $f$ contains at most one critical value, 
  and thus the conclusion follows from case~\ref{item:lctwo_onecrit}.
\end{proof}

\section{Examples with non-locally connected Julia sets}\label{sec:counterexamples}
\begin{proof}[Verification of the properties of Example \ref{ex:cossqrt}]
  It is elementary to check that $f$ has no asymptotic values and the three stated critical values, and that all critical points of $f$ are real. 
   (This also follows
    from the more general discussion that follows below.) By considering the graph of the
    restriction of $f$ to the real axis, it is easy to check that $0$ and $1$ are fixed, and that there is a unique 
    repelling fixed point $p_u$ in the interval $(0,1)$ (see Figure \ref{fig:cossqrt}). In particular, the immediate basin of $0$ contains no other
    critical points, and hence is a Jordan domain by Theorem \ref{thm:jordan}. If $u$ is chosen such that 
    $c_u>p_u$, then  
         the immediate basin of $1$ contains all positive critical points, and is hence unbounded,
         and its boundary is not locally connected at any finite point by Theorem \ref{thm:bakerweinreich}. It can be verified numerically that this is 
         the case for $u> 2.7981186\dots$\,.
\end{proof}

\subsection*{The Maclane-Vinberg method}
We shall now construct Examples \ref{example1} and \ref{example2}, of certain
  hyperbolic functions with non-locally connected Julia sets. We will use a general
method to construct real entire functions with a preassigned
sequence of critical values. We follow Eremenko and Sodin \cite{eremenkosodin1992}
in the description of the method. They credit Maclane \cite{maclane1947}
for this method and Vinberg \cite{vinberg1989} for a modern exposition thereof.  
  Another discussion of the construction can be found in \cite{eremenkoyuditskii}.
  

Let $\underline{c}=(c_n)_{n\in\Z}$ be a sequence 
satisfying $(-1)^n c_n\geq 0$ for all $n\in\Z$ and let 
\[
\Omega= \Omega(\underline{c})=
\C\backslash \bigcup_{n\in\Z} \{x+in\pi\colon -\infty< x\leq \log|c_n|\},
\]
where we set $ \{x+in\pi\colon -\infty< x\leq \log|c_n|\}=\emptyset$ if $c_n=0$.
We assume that not all $c_n$ are equal to $0$, so that $\Omega\neq \C$.
Then there exists a conformal map $\varphi$ mapping the lower half-plane
${\mathbb H}^-=\{z\in\C\colon \im z< 0\}$ onto $\Omega(\underline{c})$ such that
$\re \varphi(iy)\to +\infty$ as $y\to-\infty$. Since $\partial \Omega(\underline{c})$ is locally connected, the map $\varphi$ extends continuously to $\R$ by the
  Carath\'eodory-Torhorst Theorem \cite[Theorem 2.1]{pommerenkeboundary}; we denote this
  extension also by $\varphi$. 
  The 
  real axis then corresponds to the slits 
  $\{x+in\pi\colon -\infty< x\leq \log|c_n|\}$ under the map $\varphi$. 
 As these slits are mapped onto the real axis by the exponential function,
we deduce from the Schwarz Reflection Principle \cite[Chapter 6]{ahlforscomplexanalysis}
 that the map $g$ given by
$g(z)=\exp \varphi(z)$ extends to an entire function. 

Note that $\varphi$ and $g$ are not uniquely determined by $\underline{c}$, 
as precomposing with a map $z\mapsto az+b$ where $a,b\in\R$ and $a>0$ 
leads to a function with the same properties. If 
 $c_0=0$, $c_{\pm1} = -1$ and $c_n=c_{-n}$ for all $n$ (which will be the case in our examples),
  we can choose $\phi$ such that 
  $\{it\colon t<0\}$ is mapped onto $\R$ and 
  $\varphi(\pm 1)=\pm i\pi$. 
With this normalisation, $g=g^{\underline{c}}$ is uniquely
   determined by $\underline{c}$. A key observation is that
   $\varphi$, and hence $g$, depend continuously on the sequence $\underline{c}$, with respect to the 
   product topology on sequences.

 It turns out that the functions $g$ obtained this way belong to the
   \emph{Laguerre-P\'olya class} $\LP$,
   which consists of all entire functions that are locally uniform limits of 
   real polynomials with only real zeros. Conversely, all functions in the class $\LP$ can be obtained
   by this procedure (we shall not use this fact). Hence we shall refer to the function
   $g$ as a \emph{Laguerre-P\'olya function} for the sequence $\underline{c}$.
   We refer to \cite[\S II.9]{obreschkoff1963} for a discussion of the class $\LP$.
   A number of arguments in our proofs could be carried out by using general results for
   the Laguerre-P\'olya class, but we prefer to 
   argue directly from the definition of $g$.

\subsection*{Initial observations and examples}
If there exists $N\in\N$ such that $c_n=0$ for $n\geq N$, then $g(x)\to 0$ 
as $x\to\infty$. Similarly, $g(x)\to 0$ as $x\to-\infty$ if 
there exists $N\in\N$ such that  $c_n=0$ for $n\leq -N$.
A simple example is $\varphi(z)=-z^2$ and $g(z)=\exp(-z^2)$ which 
corresponds to $c_0=1$ and $c_n=0$ for all $n\neq 0$. If $f$ is the function
from Example~\ref{ex:cossqrt}, then $1-f$ is a Laguerre-P\'olya function, corresponding to 
$c_0=1$, $c_n=0$ if $n$ is odd and $c_n=2/(u+1)$ for $n$ even and nonzero.

Another example, which will recur in our proofs,
is given by $c_{\pm 1}=-1$ and $c_n=0$ for $|n|\neq 1$.
In this case the domain $\Omega(\underline{c})$ is given by 
$\Omega_0
=\C\backslash\{x\pm i\pi\colon  x\leq 0\}$. 
Denote the corresponding Laguerre-P\'olya function (normalised as above) by
   $g_0:=g^{\underline{c}}$; then 
\[
g_0(z)=\exp\varphi_0(z) =-z^2\exp(-z^2+1)
\]
is precisely the function from Figure~\ref{fig:asymptotic}\subref{subfig:maclanevinberglimit}. 

The critical values of a Laguerre-P\'olya function $g$ are 
precisely the $c_n$, except that $0$ is a critical value only if $c_l=0$ for some 
$l\in\Z$ for which there exist $k,m\in\Z$ with $k<l<m$, $c_k\neq 0$ and $c_m\neq 0$.
Moreover, there are critical points $\xi_n$ with $g(\xi_n)=c_n$ such that 
$\xi_n\leq \xi_{n+1}$ for all $n\in\Z$, and $g$ has no further critical
points (since $\phi$ and $\exp$ have no critical points). The limit $\lim_{x\to\infty} g(x)$ exists if and only if 
$\lim_{n\to\infty}c_n=0$, and in this case $\lim_{x\to\infty} g(x)=0$.
An analogous remark applies to the limit $\lim_{x\to-\infty} g(x)$.

\begin{figure}[hbt!]
{\def\svgwidth{\textwidth}
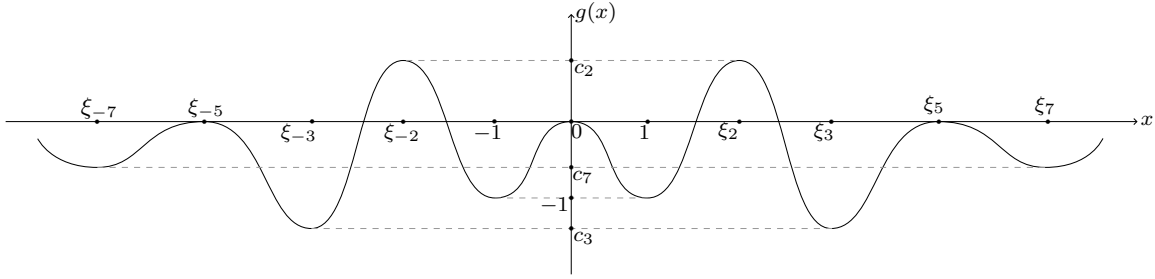}
\caption{\label{fig:sketch} \small Sketch of a Laguerre-P\'olya function $g$ on the real line, with the choices $c_0=0$, $c_{\pm 1}=-1$, $c_{\pm4}=c_{\pm5}=c_{\pm6}=0$ and $c_n=c_{-n}$ for all $n\in \N$. Consequently we may normalise so that $\xi_{\pm1}=\pm 1$ and $\xi_0=0$.  We also have multiple critical points  $\xi_{4}=\xi_{5}=\xi_{6}$ and $\xi_{-4}=\xi_{-5}=\xi_{-6}$. }
\end{figure}


We mention that the construction can be modified if $c_n$ is not defined
for all $n\in\Z$, but only for $n\leq N$ or for $n\geq M$. We can think of this as a
limit case, where $|c_{N+1}|=\infty$ or $|c_{M-1}|=\infty$, and obtain a function $g$ with $\lim_{x\to\infty} |g(x)|=\infty$ or 
$\lim_{x\to-\infty} |g(x)|=\infty$, respectively.
We shall not need these considerations.

If  $\lim_{n\to+\infty} c_n = 0$ or  $\lim_{n\to-\infty} c_n = 0$, then
  $0$ is an asymptotic value of $g$. The following shows
  that the converse also holds.

\begin{lem} \label{lemma1}
Let $g$ be a Laguerre-P\'olya function and suppose that $g$ has an asymptotic value $a\in\C$. Then $a=0$ and
\begin{equation}\label{limg}
\lim_{x\to\infty} g(x)=0
\quad\text{or}\quad
\lim_{x\to-\infty} g(x)=0.
\end{equation}
\end{lem}
\begin{proof}
 Let $a\in\C$, and let $D$ be a small disc around $a$. If $a\neq 0$, we may assume that
    $0\notin D$. Then every connected component of
    $\exp^{-1}(D)$ is bounded and intersects at most one of the lines in the complement
    of $\Omega$. Thus every component of
     \[ \phi^{-1}(\exp^{-1}(D)) = g^{-1}(D)\cap \H^- \]
     is bounded and its closure intersects the real line in at most one interval. 
    It follows (by also considering the preimage of $\overline{D}$, the reflection of $D$ in 
     the real axis) that every connected component of $g^{-1}(D)$ is bounded, and hence
    $a$ is not an asymptotic value. 

  If $a=0$, then $\exp^{-1}(D)$ is a left half-plane, and hence unbounded. However, unless
    the conclusion of the lemma is satisfied, we have
     \[ C := \min\left( \liminf_{n\to+\infty} |c_n| , \liminf_{n\to-\infty} |c_n| \right) >0. \] 
    So we can assume that the radius of $D$ was chosen smaller than $C$. 
    Then every component of
    $\exp^{-1}(D)\cap \Omega$ has bounded imaginary parts, and again
    every connected component of $g^{-1}(D)\cap\H^-$ is bounded. Since
    $g^{-1}(D)$ is symmetric with respect to the real axis, we are done.
\end{proof}

We remark that, in particular, the Maclane-Vinberg method allows us to construct uncountably many functions with two critical values that
  differ from each other in an essential manner.
\begin{obs}[Topologically inequivalent functions]\label{obs:uncountable}
  Let $A\subset \{2,4,6,\dots,\infty\}$ be nonempty. Then there exists an entire function $f\colon \C\to\C$ with
   $\sing(f^{-1})=\{0,1\}$ such that $f$ has only simple critical points over $1$, and such that $A$ is precisely the set of
   local degrees of the preimages of $0$. (Here we take $\infty\in A$ to mean that $f$ has an asymptotic value over $0$.)

  Functions corresponding to different choices of $A$ cannot be obtained from one another by pre-\ and post-composition
   with plane homeomorphisms.
\end{obs}
\begin{proof}
 Let $B\subset 2\Z$ be a set of even integers with $0\in B$ such that the length of every segment of consecutive integers in $\Z\setminus B$
  belongs to the set $A$, and such that for every element of $A$ there is a segment of this length. The desired function is obtained from the
  Maclane-Vinberg method by choosing $c_n=1$ for $n\in B$ and $c_n=0$ otherwise. The final claim follows from the fact that
  the order of a critical point is preserved under pre- and post-composition with plane homeomorphisms.
\end{proof}

\subsection*{Non-locally connected Julia sets of Laguerre-P\'olya functions} 
 We are now ready to construct the desired examples.

\begin{proof}[Construction of Example \ref{example1}]
Let $0\leq\delta<1$.
Define 
$(c_n)_{n\in\Z}$ by $c_n=0$ if $n$ is even, $c_{\pm 1}=-1$ and $c_n=-\delta$
if $n$ is odd and $|n|\geq 3$.
Let $g_\delta := g^{\overline{c}}$ be the corresponding Laguerre-P\'olya function.
Recall that $\varphi$ maps $\{it\colon t<0\}$ to $\R$
 and $\varphi(\pm 1)=\pm i\pi$.
Denoting by  $(\xi_n)_{n\in\Z}$ the sequence of critical points, with
$\xi_n\leq \xi_{n+1}$ and $g(\xi_n)=c_n$, we then have $\xi_0=0$ and
$\xi_{\pm 1}=\pm 1$ so that $g_\delta(0)=0$ and $g_\delta(\pm 1)=c_{\pm 1}=-1$.
Moreover, $g_\delta$ is an even function.
The case $\delta=0$, with $g_0(z)=-z^2\exp(-z^2+1)$,
was already considered in our description of the method.
By continuity of the Maclane-Vinberg method, we have 
\[
\lim_{\delta \to 0} g_\delta(z)=g_0(z)= -z^2 \exp(-z^2+1),
\]
locally uniformly for $z\in\C$.
Now $g_\delta(z)$ and $g_0$ have superattracting fixed
points at $0$ and $-1$. Hence we can choose $\delta$ sufficiently
small to ensure that $-\delta$ is in the immediate 
basin of~$0$ for the map $g_{\delta}$.

The function $g_{\delta}$ has no asymptotic values 
by Lemma~\ref{lemma1}, and
the only critical points are the $\xi_n$. Since $g(1)=-1$, we also see that $1$ is not in the basin of $0$. 
As the immediate attracting basins of $0$ and  $-1$ are simply-connected and
symmetric with respect to the real axis, this implies that
$\xi_0=0$ is the only critical point in the immediate basin of~$0$
and  $\xi_{-1}=-1$ is the only critical point in the immediate basin of $-1$.
Since all three critical values tend to $0$ or $-1$ under iteration, $g_{\delta}$ is hyperbolic. Hence,  by Theorem~\ref{thm:jordan}, the immediate attracting basins of $0$ and $-1$ are bounded by
Jordan curves. By assumption, $-\delta$ is contained in the immediate attracting 
basin of~$0$.  Using again that this immediate basin is simply-connected and symmetric with
respect to the real axis, we see that it actually contains the interval $[-\delta,0]$.
Now $f([\xi_2,\infty))\subset [0,\delta]$ which implies that $[\xi_2,\infty)$ is 
contained in a component of the preimage of the immediate basin of~$0$.
We conclude that $g$ satisfies the conclusion with $w_1=0$, $w_2=-1$ and $w_3=-\delta$.
\end{proof}%
\begin{proof}[Construction of Example \ref{example2}]
Recall that our goal is to construct a hyperbolic function with critical values $0$ and $1$  and no asymptotic values
  such that every Fatou component is a Jordan domain, but the Julia set is not locally connected.
   It will be slightly more convenient to normalise such that the
   critical values are $0$ and $-1$
   instead (conjugation by $z\mapsto -z$ yields the original normalisation).
   
  We begin by outlining the construction, which is based on the idea 
   that a critical point of sufficiently high degree can be
   used to approximate the behaviour of an asymptotic tract. Indeed, suppose that we
  start with the Laguerre-P\'olya function $g_0$ from the introduction to this section
  (where $c_{\pm 1}=-1$ and $c_n=0$ if $|n|\neq 1$). The super-attracting point $0$ is
  an asymptotic value for $g_0$, and since $J(g_0)$ intersects the unit disc, there 
  is an unbounded Fatou component of $g_0$ that intersects the unit circle
  $\Gamma_1 :=\{z\colon |z|=1\}$. Let us
  modify the sequence $\underline{c}$ by introducing the additional non-zero points
  $c_{\pm N}=-1$, for some large integer $N$. By continuity of the Maclane-Vinberg
   method, the corresponding function $g_1$
  is close to $g_0$, but has an additional pair of critical points of degree
   $N-1$. 
   It follows (see below for details) that $g_1$ can be chosen to have a \emph{bounded} 
  Fatou component
  that intersects both $\Gamma_1$ and $\Gamma_2 := \{z\colon |z|=2\}$. As $g_1$ still
  has an asymptotic value over $0$, we can repeat the procedure and create a function 
  $g_2$ that has two Fatou components both intersecting $\Gamma_1$ and $\Gamma_2$.
  Continuing inductively, in the limit
  we obtain a function $g$ that has no asymptotic values by Lemma \ref{lemma1},
  but has infinitely many Fatou components that intersect both $\Gamma_1$ and $\Gamma_2$.  
  Then $J(g)$ is not locally connected by Lemma~\ref{whyburnlemma}. 
\begin{figure}
 \includegraphics[width=\textwidth]{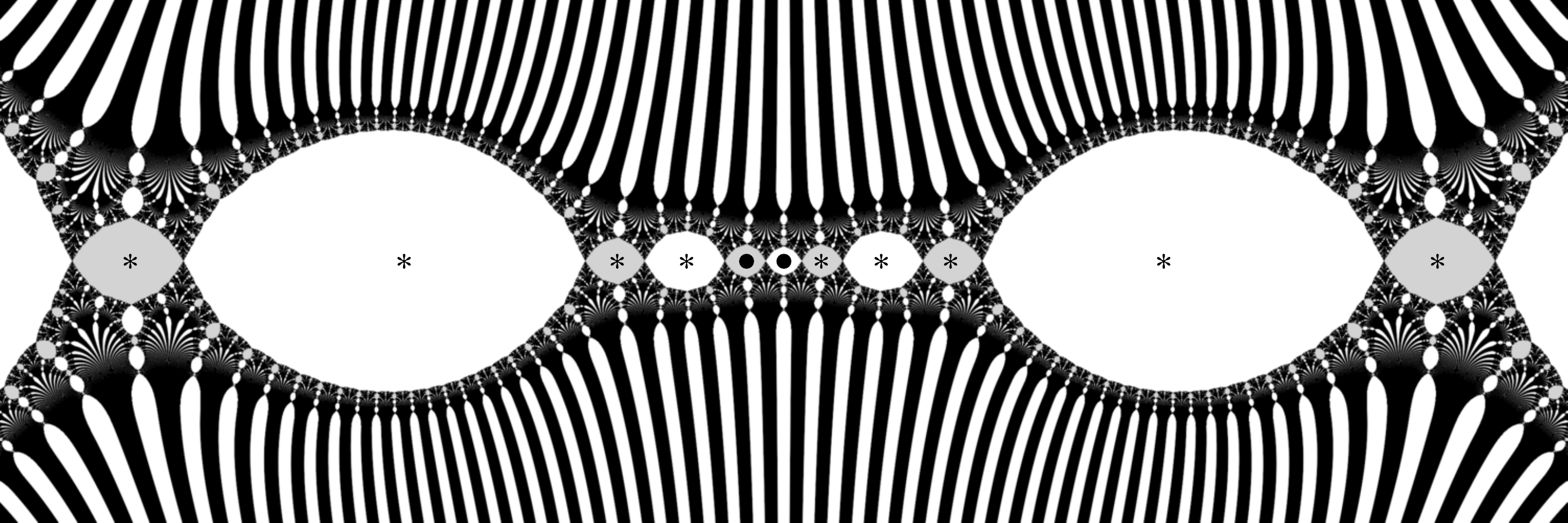}
 {\def\svgwidth{0.95\textwidth}
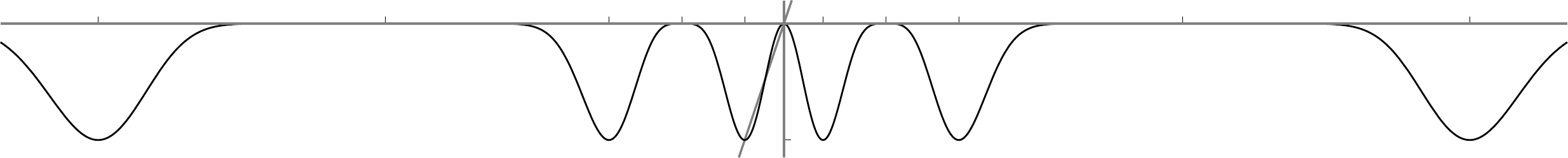}
\caption{\label{fig:maclanevinberg-5-25} \small Illustration of the construction of Example \ref{example2}. 
  Shown are the Julia set and the graph of the function $g_2$ that would arise from the choice of
   $N_1=5$ and $N_2=25$. Note the large size of the Fatou components containing 
   high-degree 
    critical points. For the actual construction of Example \ref{example2} the sequence
$(N_k)$ has to grow much more rapidly than indicated by 
the above values of $N_1$ and $N_2$.}
\end{figure}

  To provide the necessary details, let 
   $\underline{N}=(N_k)_{k\geq 0}$ be a (rapidly) increasing sequence 
     of odd positive integers with $N_0=1$.
  We define sequences $\underline{c}^K$ (depending on $\underline{N}$) by
    \[ c_n^K := \begin{cases} -1 &\text{if }|n|=N_k\text{ for some $0\leq k \leq K$} \\
                        0 &\text{otherwise},\end{cases} \]
   and their limit $\underline{c}=\underline{c}(\underline{N})$,
   \[ c_n := \begin{cases} -1 &\text{if }|n|=N_k\text{ for some $k\geq 0$} \\
                        0 &\text{otherwise}.\end{cases} \]
   Let $g_K := g^{\underline{c}^K}$ and  $g=g_{\underline{N}}:= g^{\underline{c}}$ 
    be the corresponding Laguerre-P\'olya functions. (See Figure~\ref{fig:maclanevinberg-5-25}.)

  The superattracting
    fixed points $0$ and $-1$ 
    are the only critical values of $g_K$ and of $g$.
    As in the construction of Example~\ref{example1}, we find that their immediate attracting
    basins are bounded by Jordan curves. 
    By Lemma \ref{lemma1}, $g$ has no asymptotic values, and hence 
    every Fatou component of $g$ is a bounded Jordan domain (Theorem \ref{thm:jordantwo}).

   If $(\xi_n)_{n\in\Z}$ is again the sequence of critical points of $g$, then
     $g(\xi_{\pm N_k})=-1$ and $g(\xi_n)=0$ when $n\neq |N_k|$ for all $k$. In particular,
     $\xi_{N_k+1}=\dots = \xi_{N_{k+1}-1}$ for all $k$. 
   Let $0=\eta_0 < \eta_1 < \eta_2 < \dots$ be the sequence
    of non-negative preimages of $0$; that is $\eta_k = \xi_{N_k-1}$. Also let 
     $U(\eta_k)$ be the Fatou component of $g$ containing $\eta_k$. Then
      $g(U(\eta_k))=U(0)$ for $k>0$. 
     Our goal is to show that we can choose the sequence $\underline{N}$ (inductively)
     so that there is a sequence $(\alpha_k)$ of Jordan arcs connecting
     $\Gamma_1$ and $\Gamma_2$ and a sequence $(m_k)$ of positive integers such that 
$g^{m_k}(\alpha_k)\subset U(\eta_k)$ and hence $g^{m_k+1}(\alpha_k)\subset U(0)$ for all~$k$.
As all $\alpha_k$ are in different Fatou components, this will complete the proof.

In order to define the above sequences, let $\eta_{K,k}$, for $0\leq k \leq K$, denote the critical point
of $g_K$ that corresponds to $\eta_k$, and let $U_K(\eta_{K,k})$ be the component of $F(g_K)$ that contains $\eta_{K,k}$.
We shall construct $\underline{N}$, $(\alpha_k)$ and $(m_k)$ inductively such that 
$g_K^{m_k}(\alpha_k)\subset U(\eta_{K,k})$ for $K\geq k$.
The construction will be such that we also have $g^{m_k}(\alpha_k)\subset U(\eta_{k})$ for all~$k$.

    Suppose that $N_0,\dots,N_{K},\alpha_1,\dots,\alpha_K$ and $m_1,\dots,m_K$ have already been chosen, for some $K\geq 0$.
As we let $N_{K+1}\to\infty$,  the continuity of the Maclane-Vinberg method yields that
$g_{L}\to g_K$ for $L\geq K+1$ and $g\to g_K$, regardless of the choices of $N_l$ for $l>K+1$.
(Here the convergence $g_{L}\to g_K$  as $N_{K+1}\to\infty$ is uniformly in~$L$.)
 Hence, by choosing $N_{K+1}$ large, we can achieve that 
$g_L^{m_k}(\alpha_k)\subset U(\eta_{L,k})$ for $L\geq K+1$, 
as well as $g^{m_k}(\alpha_k)\subset U(\eta_{k})$, for $k=1,\dots,K$.
      Recall that $g_K(x)\to 0$ as $x\to \infty$, so there is
     $X>0$ such that $[X,\infty)$ is contained in the basin of attraction of $0$ (for $g_K$). 
     Since $0$ and $-1$ are superattracting fixed points of $g_K$, the Julia set $J(g_K)$ 
     intersects the unit disc $\D$. It follows that there is $m_{K+1}>0$ and 
     a connected component 
     of $g_K^{-m_{K+1}}\bigl([X,\infty)\bigr)$ 
     that connects a point in $\D$ to $\infty$. Let $\alpha_{K+1}$ be a piece of
     this curve that connects $\Gamma_1$ to~$\Gamma_2$. 

Since $g_{L}\to g_K$ as $N_{K+1}\to\infty$, we have
$\eta_{L,K+1}\to\infty$ as $N_{K+1}\to\infty$, uniformly in~$L$, 
as well as $\eta_{K+1}\to\infty$.
Hence, if $N_{K+1}$ is chosen sufficiently large, then 
     $g_{L}^{m_{K+1}}(\alpha_{K+1})\subset U(\eta_{L,K+1})$ and
     $g^{m_{K+1}}(\alpha_{K+1})\subset U(\eta_{K+1})$. 
This completes the inductive construction of Example~\ref{example2}.
\end{proof}

In both Example \ref{example1} and Example \ref{example2}, 
   we constructed a function having two superattracting cycles, at 
    $0$ and at $-1$. Recall that, in both cases, local connectivity of the Julia set
    failed due to the
    preimage components of the immediate basin of $0$. The role of the fixed point at $-1$, 
    and its preimages, was to separate $0$ from all its preimages, and hence ensure that
    the immediate basins of attraction are bounded. 
 
  We remark that it is possible to modify the constructions to create a map having only 
    a single superattracting fixed point. This is achieved by normalizing our maps $g$
    so that $\pm 1$ are mapped not to $-1$, but to the first negative preimage of $0$. 
    This ensures that $g$ has a repelling fixed point between $0$ and $-1$, and the
    remainder of the proofs goes through as before.

\bibliographystyle{amsalpha}
\bibliography{biblio_lc_firstnames}

\def\cprime{$'$}
\providecommand{\bysame}{\leavevmode\hbox to3em{\hrulefill}\thinspace}
\providecommand{\MR}{\relax\ifhmode\unskip\space\fi MR }
\providecommand{\MRhref}[2]{%
  \href{http://www.ams.org/mathscinet-getitem?mr=#1}{#2}
}
\providecommand{\href}[2]{#2}
\begin{thebibliography}{RRRS11}

\bibitem[Ahl78]{ahlforscomplexanalysis}
Lars~V. Ahlfors, \emph{Complex analysis}, third ed., McGraw-Hill Book Co., New
  York, 1978, An introduction to the theory of analytic functions of one
  complex variable, International Series in Pure and Applied Mathematics.

\bibitem[AO93]{aartsoversteegen}
Jan~M. Aarts and Lex~G. Oversteegen, \emph{The geometry of {J}ulia sets},
  Trans. Amer. Math. Soc. \textbf{338} (1993), no.~2, 897--918.

\bibitem[Bak70]{bak70}
I.~Noel Baker, \emph{Limit functions and sets of non-normality in iteration
  theory}, Ann. Acad. Sci. Fenn. Ser. A I Math. \textbf{467} (1970), 11.

\bibitem[BD99]{BakDom99}
I.~Noel Baker and Patricia Dom{\'{\i}}nguez, \emph{Boundaries of unbounded
  {F}atou components of entire functions}, Ann. Acad. Sci. Fenn. Math.
  \textbf{24} (1999), no.~2, 437--464.

\bibitem[BD00a]{bakdom2000}
\bysame, \emph{Some connectedness properties of {J}ulia sets}, Complex
  Variables Theory Appl. \textbf{41} (2000), no.~4, 371--389.

\bibitem[BD00b]{bhattacharjee}
Ranjit Bhattacharjee and Robert~L. Devaney, \emph{Tying hairs for structurally
  stable exponentials}, Ergodic Theory Dynam. Systems \textbf{20} (2000),
  no.~6, 1603--1617.

\bibitem[Ber93]{waltersurvey}
Walter Bergweiler, \emph{Iteration of meromorphic functions}, Bull. Amer. Math.
  Soc. (N.S.) \textbf{29} (1993), no.~2, 151--188.

\bibitem[Bis15]{bishopfolding}
Christopher~J. Bishop, \emph{Constructing entire functions by quasiconformal
  folding}, Acta Math. \textbf{215} (2015), no.~1, 1--60.

\bibitem[BM02]{bergweilermorosawa}
Walter Bergweiler and Shunsuke Morosawa, \emph{Semihyperbolic entire
  functions}, Nonlinearity \textbf{15} (2002), no.~5, 1673--1684.

\bibitem[BM07]{beardonminda}
Alan~F. Beardon and David Minda, \emph{The hyperbolic metric and geometric
  function theory}, Quasiconformal mappings and their applications, Narosa, New
  Delhi, 2007, pp.~9--56.

\bibitem[BR84]{BakerRippon84}
I.~Noel Baker and Philip~J. Rippon, \emph{Iteration of exponential functions},
  Ann. Acad. Sci. Fenn. Ser. A I Math. \textbf{9} (1984), 49--77.

\bibitem[BW91]{bakerweinreich}
I.~Noel Baker and Jonathan~M. Weinreich, \emph{Boundaries which arise in the
  dynamics of entire functions}, Rev. Roumaine Math. Pures Appl. \textbf{36}
  (1991), no.~7-8, 413--420, Analyse complexe (Bucharest, 1989).

\bibitem[CP03]{carmonapommerenke}
Joan~J. Carmona and Christian Pommerenke, \emph{Decomposition of continua and
  prime ends}, Comput. Methods Funct. Theory \textbf{3} (2003), no.~1-2,
  253--272.

\bibitem[Den13]{deniz2013}
Asli Deniz, \emph{Entire trancendental maps with two singular values}, Ph.D.
  thesis, Roskilde University and Universitat de Barcelona, 2013.

\bibitem[Dev84]{Devaney84}
Robert~L. Devaney, \emph{Julia sets and bifurcation diagrams for exponential
  maps}, Bull. Amer. Math. Soc. (N.S.) \textbf{11} (1984), no.~1, 167--171.

\bibitem[DF08]{domfag}
Patricia Dom{\'{\i}}nguez and N{\'u}ria Fagella, \emph{Residual {J}ulia sets of
  rational and transcendental functions}, Transcendental dynamics and complex
  analysis, London Math. Soc. Lecture Note Ser., vol. 348, Cambridge Univ.
  Press, Cambridge, 2008, pp.~138--164.

\bibitem[DH85]{douhub85}
Adrien Douady and John~Hamal Hubbard, \emph{On the dynamics of polynomial-like
  mappings}, Ann. Sci. \'Ecole Norm. Sup. (4) \textbf{18} (1985), no.~2,
  287--343.

\bibitem[Dou93]{douadypincheddisk}
Adrien Douady, \emph{Descriptions of compact sets in {${\bf C}$}}, Topological
  methods in modern mathematics ({S}tony {B}rook, {NY}, 1991), Publish or
  Perish, Houston, TX, 1993, pp.~429--465.

\bibitem[Dra07]{Drasin07}
David Drasin, \emph{Regularity of growth and the class {$\mathcal{S}$}},
  Conform. Geom. Dyn. \textbf{11} (2007), 90--100.

\bibitem[DS02]{dominguezsienra}
Patricia Dom{\'{\i}}nguez and Guillermo Sienra, \emph{A study of the dynamics
  of {$\lambda\sin z$}}, Internat. J. Bifur. Chaos Appl. Sci. Engrg.
  \textbf{12} (2002), no.~12, 2869--2883.

\bibitem[EL87]{elexamples}
Alexandre~E. Eremenko and Mikhail~Yu. Lyubich, \emph{Examples of entire
  functions with pathological dynamics}, J. London Math. Soc. (2) \textbf{36}
  (1987), no.~3, 458--468.

\bibitem[EL92]{eremenko_lyubich_2}
\bysame, \emph{Dynamical properties of some classes of entire functions}, Ann.
  Inst. Fourier Grenoble \textbf{42} (1992), no.~4, 989--1020.

\bibitem[ES92]{eremenkosodin1992}
Alexandre~E. Eremenko and Mikhail~L. Sodin, \emph{Parametrization of entire
  functions of sine-type by their critical values}, Entire and subharmonic
  functions, Adv. Soviet Math., vol.~11, Amer. Math. Soc., Providence, RI,
  1992, pp.~237--242.

\bibitem[EY12]{eremenkoyuditskii}
Alexandre Eremenko and Peter Yuditskii, \emph{Comb functions}, Recent advances
  in orthogonal polynomials, special functions, and their applications,
  Contemp. Math., vol. 578, Amer. Math. Soc., Providence, RI, 2012,
  pp.~99--118.

\bibitem[Fag95]{nurialimitingstandard}
N{\'u}ria Fagella, \emph{Limiting dynamics for the complex standard family},
  Internat. J. Bifur. Chaos Appl. Sci. Engrg. \textbf{5} (1995), no.~3,
  673--699.

\bibitem[Fat20]{fatoumemoir2}
Pierre Fatou, \emph{{Sur les \'equations fonctionnelles, II.}}, {Bull. Soc.
  Math. Fr.} \textbf{48} (1920), 33--94 (French).

\bibitem[GO08]{goldbergostrovskii}
Anatoly~A. Goldberg and Iossif~V. Ostrovskii, \emph{Value distribution of
  meromorphic functions}, Translations of Mathematical Monographs, vol. 236,
  American Mathematical Society, Providence, RI, 2008, Translated from the 1970
  Russian original by Mikhail Ostrovskii, With an appendix by Alexandre
  Eremenko and James K. Langley.

\bibitem[G{\'S}97]{graczykswiatek}
Jacek Graczyk and Grzegorz {\'S}wiatek, \emph{Generic hyperbolicity in the
  logistic family}, Ann. of Math. (2) \textbf{146} (1997), no.~1, 1--52.

\bibitem[{Hei}57]{heins}
Maurice {Heins}, \emph{{Asymptotic spots of entire and meromorphic
  functions.}}, {Ann. of Math. (2)} \textbf{66} (1957), 430--439 (English).

\bibitem[KSvS07]{KSS2}
Oleg Kozlovski, Weixiao Shen, and Sebastian van Strien, \emph{Density of
  hyperbolicity in dimension one}, Ann. of Math. (2) \textbf{166} (2007),
  no.~1, 145--182.

\bibitem[Lyu97]{lyubichdensity}
Mikhail Lyubich, \emph{Dynamics of quadratic polynomials. {I}, {II}}, Acta
  Math. \textbf{178} (1997), no.~2, 185--247, 247--297.

\bibitem[Mac47]{maclane1947}
Gerald~R. MacLane, \emph{Concerning the uniformization of certain {R}iemann
  surfaces allied to the inverse-cosine and inverse-gamma surfaces}, Trans.
  Amer. Math. Soc. \textbf{62} (1947), 99--113.

\bibitem[MB10]{helenalanding}
Helena Mihaljevi{\'c}-Brandt, \emph{A landing theorem for dynamic rays of
  geometrically finite entire functions}, J. Lond. Math. Soc. (2) \textbf{81}
  (2010), no.~3, 696--714.

\bibitem[MB12]{helenaorbifolds}
\bysame, \emph{Semiconjugacies, pinched {C}antor bouquets and hyperbolic
  orbifolds}, Trans. Amer. Math. Soc. \textbf{364} (2012), no.~8, 4053--4083.

\bibitem[Mer08]{Merenkov08}
Sergei Merenkov, \emph{Rapidly growing entire functions with three singular
  values}, Illinois J. Math. \textbf{52} (2008), no.~2, 473--491.

\bibitem[Mor99]{morosawa}
Shunsuke Morosawa, \emph{Local connectedness of {J}ulia sets for transcendental
  entire functions}, Nonlinear analysis and convex analysis ({N}iigata, 1998),
  World Sci. Publ., River Edge, NJ, 1999, pp.~266--273.

\bibitem[MU08]{mayerurbanski}
Volker Mayer and Mariusz Urba{\'n}ski, \emph{Geometric thermodynamic formalism
  and real analyticity for meromorphic functions of finite order}, Ergodic
  Theory Dynam. Systems \textbf{28} (2008), no.~3, 915--946.

\bibitem[Nad92]{nadler}
Sam~B. Nadler, Jr., \emph{Continuum theory}, Monographs and Textbooks in Pure
  and Applied Mathematics, vol. 158, Marcel Dekker Inc., New York, 1992, An
  introduction.

\bibitem[Obr63]{obreschkoff1963}
Nikola Obreschkoff, \emph{Verteilung und {B}erechnung der {N}ullstellen reeller
  {P}olynome}, VEB Deutscher Verlag der Wissenschaften, Berlin, 1963.

\bibitem[Osb13]{osborne}
John~W. Osborne, \emph{Spiders' webs and locally connected {J}ulia sets of
  transcendental entire functions}, Ergodic Theory Dynam. Systems \textbf{33}
  (2013), no.~4, 1146--1161.

\bibitem[Pom92]{pommerenkeboundary}
Christian Pommerenke, \emph{Boundary behaviour of conformal maps}, Grundlehren
  der Mathematischen Wissenschaften [Fundamental Principles of Mathematical
  Sciences], vol. 299, Springer-Verlag, Berlin, 1992.

\bibitem[Rem06]{topescaping}
Lasse Rempe, \emph{Topological dynamics of exponential maps on their escaping
  sets}, Ergodic Theory Dynam. Systems \textbf{26} (2006), no.~6, 1939--1975.

\bibitem[Rem08]{localconnectivity}
\bysame, \emph{On prime ends and local connectivity}, Bull. Lond. Math. Soc.
  \textbf{40} (2008), no.~5, 817--826, Updated version available at
  http://arxiv.org/abs/math/0309022.

\bibitem[Rem09]{rigidity}
\bysame, \emph{Rigidity of escaping dynamics for transcendental entire
  functions}, Acta Math. \textbf{203} (2009), no.~2, 235--267.

\bibitem[Rem13]{fullhypdimension}
Lasse Rempe{-Gillen}, \emph{Hyperbolic entire functions with full hyperbolic
  dimension and approximation by {E}remenko--{L}yubich functions}, Proc. London
  Math. Soc. (2013).

\bibitem[Rem14]{arclike}
\bysame, \emph{Arc-like continua, {J}ulia sets of entire functions, and
  {E}remenko's {C}onjecture}, Preprint, 2014.

\bibitem[RRRS11]{strahlen}
G{\"u}nter Rottenfu{\ss}er, Johannes R{\"u}ckert, Lasse Rempe, and Dierk
  Schleicher, \emph{Dynamic rays of bounded-type entire functions}, Ann. of
  Math. (2) \textbf{173} (2011), no.~1, 77--125.

\bibitem[RS99]{ripponstallardhyperbolic}
Philip~J. Rippon and Gwyneth~M. Stallard, \emph{Iteration of a class of
  hyperbolic meromorphic functions}, Proc. Amer. Math. Soc. \textbf{127}
  (1999), no.~11, 3251--3258.

\bibitem[RS09]{RemSch09}
Lasse Rempe and Dierk Schleicher, \emph{Bifurcations in the space of
  exponential maps}, Invent. Math. \textbf{175} (2009), no.~1, 103--135.

\bibitem[RvS15]{RGvS}
Lasse Rempe{-Gillen} and Sebastian van Strien, \emph{Density of hyperbolicity
  for classes of real transcendental entire functions and circle maps}, Duke
  Math. J. \textbf{164} (2015), no.~6, 1079--1137.

\bibitem[RY08]{roeschyin}
Pascale Roesch and Yongcheng Yin, \emph{The boundary of bounded polynomial
  {F}atou components}, C. R. Math. Acad. Sci. Paris \textbf{346} (2008),
  no.~15-16, 877--880.

\bibitem[Ste93]{steinmetz}
Norbert Steinmetz, \emph{Rational iteration}, de Gruyter Studies in
  Mathematics, vol.~16, Walter de Gruyter \& Co., Berlin, 1993, Complex
  analytic dynamical systems.

\bibitem[SZ03]{schleicherzimmerperiodic}
Dierk Schleicher and Johannes Zimmer, \emph{Periodic points and dynamic rays of
  exponential maps}, Ann. Acad. Sci. Fenn. Math. \textbf{28} (2003), no.~2,
  327--354.

\bibitem[Vin89]{vinberg1989}
{\`E}rnest~B. Vinberg, \emph{Real entire functions with prescribed critical
  values}, Problems in group theory and in homological algebra ({R}ussian),
  Yaroslav. Gos. Univ., Yaroslavl\cprime, 1989, pp.~127--138.

\bibitem[Why42]{whyburn1942}
Gordon~Thomas Whyburn, \emph{Analytic {T}opology}, American Mathematical
  Society Colloquium Publications, v. 28, American Mathematical Society, New
  York, 1942.

\bibitem[Zha09]{zhangsinecharacterization}
Gaofei Zhang, \emph{Topological characterization of the hyperbolic maps in the
  sine family}, Preprint arxiv:0904.4081, 2009.

\end{thebibliography}
\end{document}